\documentclass[11pt,a4paper]{amsart}
\usepackage[T1]{fontenc}
\usepackage{lmodern}
\usepackage{amsmath,amssymb,amsthm,mathtools}
\usepackage[margin=28mm]{geometry}
\usepackage[hidelinks]{hyperref}
\hypersetup{pdftitle={Sharp Holder regularity for autonomous Beltrami equations},
	pdfauthor={Jie Huang and Jian-Feng Zhu},
	pdfkeywords={autonomous Beltrami equation, sharp Holder regularity, critical flux}}
\numberwithin{equation}{section}
\newtheorem{theorem}{Theorem}[section]
\newtheorem{proposition}[theorem]{Proposition}
\newtheorem{lemma}[theorem]{Lemma}
\newtheorem{corollary}[theorem]{Corollary}
\newtheorem*{theoremA}{Theorem A}
\newtheorem*{conjecture}{Conjecture}
\theoremstyle{remark}

\newcommand{\C}{\mathbb C}
\newcommand{\R}{\mathbb R}
\newcommand{\Z}{\mathbb Z}
\newcommand{\T}{\mathbb R/(2\pi\mathbb Z)}
\newcommand{\loc}{\mathrm{loc}}
\newcommand{\dd}{\mathrm{d}}
\newcommand{\dA}{\,\mathrm{d}A}
\newcommand{\HH}{\mathcal H}
\newcommand{\FF}{\mathcal F}
\newcommand{\DD}{\mathcal D}
\newcommand{\BB}{\mathcal B}
\newcommand{\MM}{\mathcal M}
\newcommand{\WW}{\mathcal W}
\DeclareMathOperator{\Rea}{Re}
\DeclareMathOperator{\Lip}{Lip}
\DeclareMathOperator{\diver}{div}

\makeatletter
\@namedef{subjclassname@2020}{\textup{2020} Mathematics Subject Classification}
\makeatother
\title[Sharp regularity for autonomous Beltrami equations]{Sharp H\"older regularity for autonomous Beltrami equations}
\author{Jie Huang}
\address{Jie Huang, School of Mathematical Sciences, Huaqiao University, Quanzhou 362021, P. R. China}
\email{jie.huang@hqu.edu.cn}
\author{Jian-Feng Zhu}
\address{Jian-Feng Zhu, Department of Mathematics, Shantou University, Shantou, Guangdong 515063, P. R. China}
\email{flandy@stu.edu.cn}
\date{October 18, 2026}
\subjclass[2020]{35J60, 30C62, 35B65}
\keywords{Autonomous Beltrami equation, sharp H\"older regularity, quasiregular mapping, Fourier--Mellin analysis, critical flux}
\begin{document}
	\begin{abstract}
Let $\Omega\subset\mathbb C$ be a domain and let $\mathcal H:\mathbb C\to\mathbb C$ be Lipschitz with constant $0<k<1$. We prove that every weak $W^{1,2}_{\mathrm{loc}}$ solution of the autonomous Beltrami equation $f_{\bar z}=\mathcal H(f_z)$ belongs to $C^{1,\gamma}_{\mathrm{loc}}$, where $\gamma=3/(2K+1)$ and $K=(1+k)/(1-k)$. The exponent is optimal, resolving the conjecture of Astala, Clop, Faraco, J\"a\"askel\"ainen and Koski. 
	\end{abstract}
	\maketitle
	
	\section{Introduction}\label{sec:introduction}
	
	Quasiregular mappings provide a geometric approach to the regularity of
	elliptic equations in the plane. A distortion inequality for the
	differential of a map gives decay of its Dirichlet energy and H\"older
	continuity. This connection
	underlies the work of Morrey \cite{Morrey} and Finn and Serrin
	\cite{FinnSerrin}, the theory of generalized analytic functions
	\cite{Vekua}, and the later development of quasiconformal methods
	\cite{AIM}. For background on regularity theory for nonlinear elliptic systems, see \cite{Giaquinta}.
	
	The relation is particularly direct in two dimensions:
	under ellipticity assumptions, a solution of a divergence-form equation
	and its local stream function combine into a map satisfying a Beltrami
	system. This correspondence allows estimates for elliptic partial
	differential equations (PDEs) to be obtained from the distortion of
	planar mappings. Quasiconformal mappings are the homeomorphic members
	of the quasiregular class, so the same estimates also describe the local
	geometry of changes of variables.
	
	H\"older continuity provides a quantitative modulus of continuity. It
	gives a continuous representative of a weak solution and bounds its
	oscillation on shrinking disks. For
	quasiconformal mappings, it controls how small distances change under
	the map. Applied to the first derivatives of a PDE solution, it gives
	$C^{1,\alpha}$ regularity and quantitative control of the gradient.
	The optimal exponent measures the strongest such control allowed by the
	ellipticity. When several derivatives of one function are quasiregular,
	their compatibility can give better regularity than the distortion
	inequality for each derivative separately.
	
	Let $\Omega\subset\C$ be a domain. We consider the autonomous nonlinear
	Beltrami equation
	\begin{equation}\label{eq:autonomous}
		f_{\bar z}=\HH(f_z)\qquad\text{a.e. in }\Omega,
	\end{equation}
	where $\HH:\C\to\C$ satisfies
	\begin{equation}\label{eq:lipschitz}
		|\HH(\zeta)-\HH(\eta)|\le k|\zeta-\eta|
		\qquad(\zeta,\eta\in\C),\qquad 0<k<1.
	\end{equation}
	By a natural weak solution we mean a map in
	$W^{1,2}_{\loc}(\Omega,\C)$ satisfying \eqref{eq:autonomous} almost
	everywhere. We use the parameters
	\begin{equation}\label{eq:parameters}
		K=\frac{1+k}{1-k},\qquad
		\gamma=\frac{3(1-k)}{3+k}=\frac{3}{2K+1}.
	\end{equation}
	The usual normalization $\HH(0)=0$ is not needed for local regularity.
	Indeed, subtracting $\HH(0)\bar z$ from $f$ replaces the structure field by
	$\HH-\HH(0)$ without changing its Lipschitz constant.
	
	Difference quotients of a solution of \eqref{eq:autonomous} satisfy the
	$K$-quasiregular differential inequality. It follows that
	$f\in W^{2,2}_{\loc}$ and that each real directional derivative of $f$
	satisfies the same inequality. See
\cite[Proposition~2.1]{ACFJK17} and \cite[Proposition~8]{ACFJK20}. The classical regularity
	theorem for quasiregular maps then gives \mbox{$f\in C^{1,1/K}_{\loc}$}.
	This estimate treats the directional derivatives independently and does
	not use the commutation of their mixed derivatives. Nonlinear Beltrami
	equations also arise in the theory of normalized quasiconformal solutions
	and families of quasiconformal mappings \cite{ACFJS,ACFJ}. Here we consider
	all natural weak solutions, without an injectivity assumption.
	
	Astala, Clop, Faraco, J\"a\"askel\"ainen and Koski \cite{ACFJK20}
	used the simultaneous distortion inequalities for the directional
	derivatives to improve the classical exponent. They proved the following
	theorem.
	
	\begin{theoremA}[{\cite[Theorem 1]{ACFJK20}}]
		Every $W^{1,2}_{\loc}$ solution of \eqref{eq:autonomous} under
		\eqref{eq:lipschitz} belongs to $C^{1,\alpha_K}_{\loc}$, where
		\begin{equation}\label{eq:old-exponent}
			\alpha_K=\min\left\{\frac{4}{3K+1},\frac{K+1}{3K-1}\right\}
			>\frac1K.
		\end{equation}
	\end{theoremA}
	
	They also constructed the solution
	$f_\gamma(z)=z^2|z|^{\gamma-1}$ for a structure field with Lipschitz
	constant $k$. Since its first derivatives are homogeneous of degree
	$\gamma$, no exponent larger than $\gamma$ can hold for all such
	equations. Immediately after their Theorem~1, they proposed the following
	conjecture \cite[p.~232]{ACFJK20}.
	
	\begin{conjecture}[Astala--Clop--Faraco--J\"a\"askel\"ainen--Koski]
		Under \eqref{eq:lipschitz}, every natural weak solution of
		\eqref{eq:autonomous} belongs to $C^{1,\gamma}_{\loc}(\Omega,\C)$,
		where $\gamma=3/(2K+1)$, and this exponent is optimal.
	\end{conjecture}
	
	For $K>1$, the known exponent satisfies $\alpha_K<\gamma$, leaving a
	gap between the regularity theorem and the homogeneous example. We prove
	the conjecture, including the endpoint.
	
	\begin{theorem}\label{thm:natural}
		Let $\HH:\C\to\C$ satisfy \eqref{eq:lipschitz}, and let
		$f\in W^{1,2}_{\loc}(\Omega,\C)$ solve \eqref{eq:autonomous}. Then
		\[
		f\in W^{2,2}_{\loc}(\Omega,\C)\cap C^{1,\gamma}_{\loc}(\Omega,\C).
		\]
		The exponent $\gamma$ in \eqref{eq:parameters} is optimal. Moreover, the
		quantitative estimates \eqref{eq:main-energy} and \eqref{eq:main-holder}
		below hold with constants depending only on $k$.
	\end{theorem}
	
	We deduce Theorem~\ref{thm:natural} from a result for maps whose real
	directional derivatives are quasiregular. For a complex-valued Sobolev
	map, let $Df$ denote its real differential and let $D^2f$ denote its
	Cartesian Hessian. We equip both with the Frobenius norm. For
	$0<a<1$, set
	\[
	[Df]_{a,U}=\sup_{\substack{x,y\in U\\x\ne y}}
	\frac{|Df(x)-Df(y)|}{|x-y|^a}.
	\]
	We write $U\Subset\Omega$ when the closure $\overline U$ is a compact
	subset of $\Omega$. In particular, the condition
	$B(c,2R)\Subset\Omega$ in the next theorem means that the whole closed
	disk $\overline{B(c,2R)}$ lies in $\Omega$ and is separated from
	$\partial\Omega$ by a positive distance. The estimates below concern
	such interior disks.
	
	\begin{theorem}\label{thm:directional}
		Suppose that $f\in W^{2,2}_{\loc}(\Omega,\C)$ and that, on a common set
		of full measure,
		\begin{equation}\label{eq:directional}
			|\partial_e f_{\bar z}|\le k|\partial_e f_z|
			\qquad\text{for every }e\in\C,\quad |e|=1.
		\end{equation}
		Then $f\in C^{1,\gamma}_{\loc}(\Omega,\C)$. There is a constant $C_k>0$
		such that, for every $B(c,2R)\Subset\Omega$ and $0<r\le R$,
		\begin{equation}\label{eq:main-energy}
			\int_{B(c,r)}|D^2f|^2\dA
			\le C_k\left(\frac rR\right)^{2\gamma}
			\int_{B(c,2R)}|D^2f|^2\dA,
		\end{equation}
		and
		\begin{equation}\label{eq:main-holder}
			[Df]_{\gamma,B(c,R)}
			\le C_k R^{-\gamma}
			\left(\int_{B(c,2R)}|D^2f|^2\dA\right)^{1/2}.
		\end{equation}
	\end{theorem}
	
	Condition \eqref{eq:directional} is precisely the $K$-quasiregular
	differential inequality for every real directional derivative. It does
	not require $f_{\bar z}$ to be a globally defined function of $f_z$.
	Theorem~\ref{thm:directional} applies to a larger class than
	the solutions of \eqref{eq:autonomous}.
	
	Improved exponents from differential compatibility have also been studied
	for quasiregular gradient mappings. Baernstein and Kovalev \cite{BK}
	proved such estimates for planar equations in nondivergence form. Hou and
	Zhu \cite{HZ} obtained a further improvement by optimizing a family of
	quadratic lower bounds for the Jacobian in the Fourier comparison for
	gradient mappings. These results illustrate how differential
	compatibility can improve the exponent obtained from distortion alone.
	Here the compatibility comes from a different condition. All real
	directional derivatives of one complex-valued function satisfy the same
	quasiregular inequality.
	
	We now describe the steps that lead from this condition to the critical
	estimate. As in \cite[Lemma~11]{ACFJK20}, put
	$u=f_{zz}$, $v=f_{z\bar z}$ and $w=f_{\bar z\bar z}$. The directional
	inequalities give
	\[
	|w|+(1-k)|v|\le k|u|.
	\]
	For a fixed center $c$, let $\vartheta=\arg(z-c)$. We retain the relative
	phase of the Hessian components in the quadratic form
	\begin{equation}\label{eq:intro-Q}
		Q_c(f)=|w-(1-k)e^{2i\vartheta}v|^2-k^2|u|^2.
	\end{equation}
	The modulus inequality above gives $Q_c(f)\le0$. The factor
	$e^{2i\vartheta}$ makes the two terms inside the first square
	have the same angular factor when the Cartesian derivatives are written
	in polar coordinates. This choice leads to an explicit quadratic form
	on the logarithmic cylinder with radial coordinate $s=\log|z-c|$.
	
	Fourier series in $\vartheta$ and the Fourier transform in $s$ reduce the
	weighted integral of \eqref{eq:intro-Q} to multiplication by a symbol.
	Factoring this symbol gives a coercive estimate for the tangential
	components of the Hessian at every exponent below $\gamma$. The coefficient of the
	fourth-order radial term vanishes, so the symbol does not control the
	full Hessian of an arbitrary test function. For maps satisfying
	\eqref{eq:directional}, the pointwise distortion inequality supplies the
	missing component. To apply the estimate to weak solutions, we average
	over the center. The resulting bounded coefficients allow strong
	$W^{2,2}$ approximation without any prior weighted integrability
	assumption. The estimate then yields weighted integrability of the Hessian.
	
	At the critical exponent $\gamma$, the symbol vanishes only for angular
	mode $m=2$ and radial frequency $\xi=0$. This degeneracy prevents a
	direct passage to the limit in the subcritical estimate. The critical
	weighted integral over a disk is itself infinite for $f_\gamma$. We
	integrate by parts on the logarithmic cylinder and keep the boundary
	terms. The resulting flux is negative definite on the resonant mode.
	The positive bulk term controls the remaining modes.
	
	The sign of this flux at the inner end $s\to-\infty$ requires a separate
	argument. After subtracting the affine Taylor polynomial of $f$ at $c$, we use one
	subcritical exponent to estimate the weighted Hessian and a second,
	slightly larger exponent to estimate the Taylor remainder. Positive flux
	at any inner radius would then force backward exponential growth,
	contradicting these bounds. This proves that the flux is nonpositive.
	Its sign controls the resonant mode on each annulus, while the bulk term
	controls all other modes. An averaged trace estimate selects the outer
	radius with a bound uniform in $c$. Summing the annular estimates proves
	\eqref{eq:main-energy}, from which \eqref{eq:main-holder} follows by the
	Morrey--Campanato theorem \cite{Morrey,Campanato}.
	
	The reduction by difference quotients, the modulus inequality, and the
	example establishing sharpness come from \cite{ACFJK20}. The argument
	here adds the angular phase in \eqref{eq:intro-Q}, the factorization of
	its Fourier--Mellin symbol, and the boundary flux estimate needed at the
	critical exponent. The estimates on annuli recover the critical energy
	decay despite the divergence of the corresponding weighted disk integral.
	
	The local theorem has two consequences. The first concerns autonomous
	equations in divergence form. Their relation to planar first-order
	systems is discussed in \cite{LN,IS,AIM}. The monotone framework goes
	back to Leray and Lions \cite{LL}, and the quasiconformal geometry of
	monotone maps is studied in \cite{Kovalev}. Astala et al.\
	\cite[Theorem~5]{ACFJK20} identified a strong ellipticity condition for
	which this correspondence preserves the parameter $K$. Under that
	condition, Theorem~\ref{thm:natural} gives the endpoint of their
	Corollary~2.
	
	\begin{corollary}\label{cor:leray-lions}
		Let $\mathbf{A}:\R^2\to\R^2$, $\mathbf{A}(0)=0$, satisfy
		\begin{equation}\label{eq:strong-ellipticity}
			|\xi-\eta|^2+|\mathbf{A}(\xi)-\mathbf{A}(\eta)|^2
			\le (K+K^{-1})\langle\xi-\eta,\mathbf{A}(\xi)-\mathbf{A}(\eta)\rangle
		\end{equation}
		for all $\xi,\eta\in\R^2$, where $K>1$. Every distributional solution
		$u\in W^{1,2}_{\loc}(\Omega,\R)$ of
		\[
		\diver\mathbf{A}(\nabla u)=0
		\]
		belongs to $C^{1,\gamma}_{\loc}(\Omega,\R)$, with
		$\gamma=3/(2K+1)$. This exponent is optimal in this class.
	\end{corollary}
	
	The stream function used in the proof is constructed on interior disks,
	so no topological assumption on $\Omega$ is necessary.
	
	The second consequence characterizes global H\"older regularity of
	entire solutions by their growth at infinity. Define
	\begin{align}
		E_f(R)&=\int_{B(0,R)}|D^2f|^2\dA,\nonumber\\
		\mathcal E_\gamma(f)&=
		\left(\liminf_{R\to\infty}R^{-2\gamma}E_f(R)\right)^{1/2},
		\label{eq:growth-energy}\\
		\mathcal V_\gamma(f)&=
		\liminf_{R\to\infty}R^{-1-\gamma}\sup_{B(0,R)}|f|.
		\label{eq:growth-value}
	\end{align}
	The local energy estimate gives the following comparison between these
	growth quantities and the global H\"older seminorm of $Df$.
	
	\begin{theorem}\label{thm:entire}
		For an entire map satisfying the hypotheses of Theorem~\ref{thm:directional},
		\begin{equation}\label{eq:global-equivalence}
			C_k^{-1}[Df]_{\gamma,\C}
			\le\mathcal E_\gamma(f)
			\le\left(\limsup_{R\to\infty}R^{-2\gamma}E_f(R)\right)^{1/2}
			\le C_k[Df]_{\gamma,\C},
		\end{equation}
		where the inequalities are understood in the extended nonnegative reals.
		If $f$ is an entire solution of \eqref{eq:autonomous}, then also
		\begin{equation}\label{eq:value-equivalence}
			C_k^{-1}[Df]_{\gamma,\C}\le\mathcal V_\gamma(f)
			\le C_k[Df]_{\gamma,\C}.
		\end{equation}
		In particular, $\mathcal E_\gamma(f)=0$ implies that $f$ is affine in
		the directional class. For autonomous solutions, the condition
		$\mathcal V_\gamma(f)=0$ also implies that $f$ is affine. The powers $2\gamma$
		and $1+\gamma$ are sharp.
	\end{theorem}
	
	Theorem~\ref{thm:entire} concerns H\"older regularity and quadratic energy
	growth. The higher integrability of derivatives is a separate question
	in planar elliptic theory \cite{Astala94,AFS}. In particular, the sharper
	Sobolev estimates asked for in \cite{ACFJK20} are not a consequence of
	the H\"older theorem proved here.
	
	Section~\ref{sec:preliminaries} gives the reduction to directional
	quasiregularity and develops the quadratic form and its Fourier--Mellin
	symbol. Section~3 contains the proofs of the main results. We first
	establish the subcritical estimates in Section~\ref{subsec:subcritical}. The flux
	identity is derived in Section~\ref{subsec:flux}, and the inner boundary
	argument and Theorem~\ref{thm:directional} are proved in
	Section~\ref{subsec:endpoint}. Section~\ref{subsec:natural} proves
	Theorem~\ref{thm:natural}, the sharpness assertion, and
	Corollary~\ref{cor:leray-lions}. Section~\ref{subsec:entire} proves
	Theorem~\ref{thm:entire}.
	
	\section{Preliminaries and the weighted quadratic form}\label{sec:preliminaries}
	
	\subsection{Real and complex derivatives}
	We identify $\C$ with $\R^2$ and write $z=x+iy$ and $f=f^1+if^2$.
	The real differential of $f$ is the linear map
	$Df:\R^2\to\C$ given by
	\[
	Df[e]=e_xf_x+e_yf_y,\qquad e=(e_x,e_y)\in\R^2.
	\]
	Equivalently, its real matrix has rows $\nabla f^1$ and $\nabla f^2$.
	The Cartesian Hessian is the symmetric, real bilinear map
	$D^2f:\R^2\times\R^2\to\C$ given by
	\[
	D^2f[e,d]
	=e_xd_xf_{xx}+(e_xd_y+e_yd_x)f_{xy}+e_yd_yf_{yy}.
	\]
	Its matrix in the coordinate basis is
	\[
	\begin{pmatrix}f_{xx}&f_{xy}\\ f_{xy}&f_{yy}\end{pmatrix}.
	\]
	This complex matrix records the two real Hessian matrices $D^2f^1$
	and $D^2f^2$, with bilinearity over $\R$. For Sobolev
	maps these matrices consist of weak derivatives and are defined almost
	everywhere.
	
	We use the Wirtinger derivatives and the area measure
	\[
	f_z=\tfrac12(f_x-if_y),\qquad
	f_{\bar z}=\tfrac12(f_x+if_y),\qquad \dd A=\dd x\,\dd y.
	\]
	The corresponding differential operators are
	$\partial_z=\tfrac12(\partial_x-i\partial_y)$ and
	$\partial_{\bar z}=\tfrac12(\partial_x+i\partial_y)$. By the
	\emph{Wirtinger Hessian} we mean the array of second derivatives
	\[
	H_{\mathrm W}f
	=\begin{pmatrix}
		f_{zz}&f_{z\bar z}\\
		f_{\bar z z}&f_{\bar z\bar z}
	\end{pmatrix}
	=\begin{pmatrix}
		f_{zz}&f_{z\bar z}\\
		f_{z\bar z}&f_{\bar z\bar z}
	\end{pmatrix}.
	\]
	Here, for example, $f_{zz}=\partial_z(\partial_zf)$ and
	$f_{z\bar z}=\partial_z(\partial_{\bar z}f)$. The mixed
	derivatives commute, also in the sense of distributions. Expanding
	the operators gives
	\[
	\begin{aligned}
		f_{zz}&=\tfrac14(f_{xx}-2if_{xy}-f_{yy}),\\
		f_{z\bar z}&=\tfrac14(f_{xx}+f_{yy})=\tfrac14\Delta f,\\
		f_{\bar z\bar z}&=\tfrac14(f_{xx}+2if_{xy}-f_{yy}).
	\end{aligned}
	\]
	Conversely, the Cartesian derivatives are recovered from
	\[
	\begin{aligned}
		f_{xx}&=f_{zz}+2f_{z\bar z}+f_{\bar z\bar z},\\
		f_{xy}&=i(f_{zz}-f_{\bar z\bar z}),\\
		f_{yy}&=-f_{zz}+2f_{z\bar z}-f_{\bar z\bar z}.
	\end{aligned}
	\]
	The two Hessians contain exactly the same information.
	The Wirtinger notation does not require $f$ to be holomorphic.
	The matrix $H_{\mathrm W}f$ is symmetric under transposition because
	its mixed entries are equal. This symmetry does not imply Hermitian
	symmetry.
	For complex-valued $f$, in general,
	$f_{\bar z\bar z}\ne\overline{f_{zz}}$ and $f_{z\bar z}$ need
	not be real. If $f$ is real-valued, then
	$f_{\bar z\bar z}=\overline{f_{zz}}$ and
	$f_{z\bar z}=\Delta f/4\in\R$.
	All these identities hold almost everywhere for
	$f\in W^{2,2}_{\loc}$.
	
	If a real direction is identified with $e=e_x+ie_y\in\C$, then
	\[
	\partial_e=e_x\partial_x+e_y\partial_y
	=e\partial_z+\bar e\partial_{\bar z}.
	\]
	Thus, the directions in \eqref{eq:directional} are ordinary real unit
	vectors. The Frobenius norms of the differential and Hessian are
	\begin{align}
		|Df|^2&=|f_x|^2+|f_y|^2
		=2(|f_z|^2+|f_{\bar z}|^2),\nonumber\\
		|D^2f|^2&=|f_{xx}|^2+2|f_{xy}|^2+|f_{yy}|^2
		=4(|f_{zz}|^2+2|f_{z\bar z}|^2+|f_{\bar z\bar z}|^2).
		\label{eq:norms}
	\end{align}
	They are the sums of the corresponding squared norms for $f^1$ and
	$f^2$ and do not depend on the choice of an orthonormal coordinate
	basis. We write $C_k$ for a positive constant depending only on $k$,
	whose value may change between occurrences. Other dependence is
	indicated by subscripts. For a measurable set $U$ of finite positive
	area and an integrable function $v$, its integral average is
	\[
	|U|=\int_U\dA,\qquad
	v_U=\frac{1}{|U|}\int_Uv(z)\,\dd A(z).
	\]
	For vector-valued $v$, the integral is taken componentwise.
	
	Fix a center $c\in\C$, and put $z-c=re^{i\vartheta}$, $r>0$. At this
	point let
	\[
	R=(\cos\vartheta,\sin\vartheta),\qquad
	T=(-\sin\vartheta,\cos\vartheta).
	\]
	The notation $f_{RR}$, $f_{RT}$ and $f_{TT}$ means
	$D^2f[R,R]$, $D^2f[R,T]$ and $D^2f[T,T]$, respectively. In each
	evaluation the vectors are fixed at the point. These quantities differ
	from repeated differentiation along a varying polar vector field.
	For example, the vector field $T$ satisfies $\partial_TT=-R/r$, so
	\[
	\partial_T(\partial_Tf)=f_{TT}-r^{-1}f_R.
	\]
	Since $\partial_R=\partial_r$ and $\partial_T=r^{-1}\partial_\vartheta$
	on scalar functions, the Hessian components are
	\begin{equation}\label{eq:polar-hessian}
		f_{RR}=f_{rr},\qquad
		f_{RT}=r^{-1}f_{r\vartheta}-r^{-2}f_\vartheta,\qquad
		f_{TT}=r^{-2}f_{\vartheta\vartheta}+r^{-1}f_r.
	\end{equation}
	For instance, if $f(z)=|z-c|^2$, then $f_{\vartheta\vartheta}=0$ but
	$f_{TT}=2$. The term $f_r/r$ in \eqref{eq:polar-hessian} cannot be
	omitted. For $W^{2,2}_{\loc}$ maps these formulas hold almost
	everywhere on annuli. We write $B_R=B(c,R)$ when the center is fixed.
	
	\subsection{Energy decay and directional quasiregularity}
	The following version of the Morrey--Campanato criterion will be
	applied to $v=Df$.
	
	\begin{lemma}\label{lem:campanato}
		Let $0<a<1$ and $v\in W^{1,2}(B(c,2R),\R^N)$. Suppose that
		\[
		\int_{B(y,\rho)}|Dv|^2\dA\le M^2\rho^{2a}
		\quad\text{for }y\in B(c,R),\quad 0<\rho\le R/2.
		\]
		Then $v$ has a representative in $C^a(B(c,R))$ and
		$[v]_{a,B(c,R)}\le C_aM$.
		The same conclusion holds on $\C$ if the hypothesis holds for all
		centers and radii.
	\end{lemma}
	\begin{proof}
		Poincar\'e's inequality and the energy hypothesis first give a bound
		for the mean oscillation:
		\[
		\frac{1}{|B(y,\rho)|}
		\int_{B(y,\rho)}|v-v_{B(y,\rho)}|^2\dA
		\le C\int_{B(y,\rho)}|Dv|^2\dA
		\le CM^2\rho^{2a}.
		\]
		The cancellation of the factor $\rho^2$ in Poincar\'e's inequality
		with the area $\pi\rho^2$ is specific to the two-dimensional domain.
		Fix $y\in B(c,R)$ and $0<\rho\le R/2$, and set
		$\rho_j=2^{-j}\rho$, $B_j=B(y,\rho_j)$ and $m_j=v_{B_j}$.
		Since $|B_j|/|B_{j+1}|=4$, the Cauchy--Schwarz inequality gives
		\begin{align*}
			|m_{j+1}-m_j|
			&\le |B_{j+1}|^{-1/2}
			\left(\int_{B_{j+1}}|v-m_j|^2\dA\right)^{1/2}\\
			&\le 2\left(\frac{1}{|B_j|}
			\int_{B_j}|v-m_j|^2\dA\right)^{1/2}
			\le CM\rho_j^a.
		\end{align*}
		The geometric series converges because $a>0$. Consequently, the
		averages have a limit, denoted by $\widetilde v(y)$, and
		\[
		|\widetilde v(y)-v_{B(y,\rho)}|
		\le CM\sum_{j=0}^{\infty}(2^{-j}\rho)^a
		=\frac{CM}{1-2^{-a}}\rho^a.
		\]
		This limit does not depend on the initial radius. Indeed, for
		$\rho_{j+1}<r\le\rho_j$ the same comparison gives
		\[
		|v_{B(y,r)}-m_j|
		\le\frac{\rho_j}{r}
		\left(\frac{1}{|B_j|}\int_{B_j}|v-m_j|^2\dA\right)^{1/2}
		\le CM\rho_j^a\to 0.
		\]
		Thus, $\widetilde v(y)$ is the limit of all ball averages as their
		radii tend to zero, and it is defined at every $y\in B(c,R)$.
		
		To compare two distinct points $x,y\in B(c,R)$, first suppose that
		$d=|x-y|\le R/8$. Both $B(x,d)$ and $B(y,d)$ lie in $B(x,2d)$,
		whose radius is at most $R/2$. Applying the mean-oscillation estimate
		on this larger ball yields
		\[
		|v_{B(x,d)}-v_{B(x,2d)}|
		+|v_{B(y,d)}-v_{B(x,2d)}|
		\le CMd^a.
		\]
		For each term, write the difference as the average of
		$v-v_{B(x,2d)}$ over the smaller ball and apply the
		Cauchy--Schwarz inequality. The ratio of the two areas is $4$.
		Combining this
		bound with the preceding limit estimate at $x$ and $y$ gives
		\[
		|\widetilde v(x)-\widetilde v(y)|\le C_aM|x-y|^a.
		\]
		If $d>R/8$, divide the segment from $x$ to $y$ into
		$n=\lceil8d/R\rceil\le16$ equal parts. All division points remain
		in the convex ball $B(c,R)$, and successive points are at distance
		at most $R/8$. Summing the estimate for these pairs gives
		$C_aMn(d/n)^a\le C_aM d^a$, after changing $C_a$.
		Lebesgue differentiation identifies $\widetilde v$ with $v$ almost
		everywhere, so it is the required representative. On $\C$, the radii
		are unrestricted, and any two points can be compared directly at
		scale $d$. See also \cite{Morrey,Campanato}.
	\end{proof}
	
	We say that $w$ satisfies the $K$-quasiregular differential inequality
	if $|w_{\bar z}|\le k|w_z|$ almost everywhere. This terminology
	includes constant maps.
	
	\begin{lemma}\label{lem:caccioppoli}
		Suppose that $w\in W^{1,2}_{\loc}(U,\C)$ and
		$|w_{\bar z}|\le k|w_z|$ almost everywhere. Then, for every real
		$\eta\in C_c^\infty(U)$ and every $q\in\C$, we have
		\begin{equation}\label{eq:caccioppoli}
			\int_U\eta^2|Dw|^2\dA
			\le\frac{2(1+k^2)}{(1-k)^2}
			\int_U|\nabla\eta|^2|w-q|^2\dA.
		\end{equation}
	\end{lemma}
	\begin{proof}
		For every compactly supported $v\in W^{1,2}(\C,\C)$, one has
		$\|v_z\|_2=\|v_{\bar z}\|_2$. In fact, for smooth maps this follows by
		integrating the Jacobian $\det Dv$, whose integral vanishes by
		integration by parts. The Sobolev case follows by density.
		Apply this identity to $v=\eta(w-q)$. Since $\eta$ is real,
		\begin{align*}
			\|\eta w_z\|_2
			&\le\|v_z\|_2+\|\eta_z(w-q)\|_2\\
			&=\|v_{\bar z}\|_2+\|\eta_z(w-q)\|_2\\
			&\le k\|\eta w_z\|_2+2\|\eta_z(w-q)\|_2.
		\end{align*}
		The relations $|\eta_z|=|\nabla\eta|/2$ and
		$|Dw|^2\le2(1+k^2)|w_z|^2$ now give the claimed estimate.
	\end{proof}
	
	\begin{lemma}\label{lem:natural}
		Every weak solution $f\in W^{1,2}_{\loc}(\Omega,\C)$ of
		\eqref{eq:autonomous} under \eqref{eq:lipschitz} belongs to
		$W^{2,2}_{\loc}$ and satisfies \eqref{eq:directional} on a common set
		of full measure.
	\end{lemma}
	\begin{proof}
		Choose open sets $U\Subset U'\Subset\Omega$ and a number $\delta>0$
		such that the $\delta$-neighborhood $U''=U'+B(0,\delta)$ has compact
		closure in $\Omega$. Here the sum denotes all points $z+w$ with
		$z\in U'$ and $|w|<\delta$. Fix a real unit direction $e$ and,
		for $0<|t|<\delta$, define the difference quotient on $U'$ by
		\[
		F_t^e(z)=\frac{f(z+te)-f(z)}{t}.
		\]
		Translations commute with weak derivatives. The equation at $z+te$
		and $z$, followed by \eqref{eq:lipschitz}, gives
		\begin{align*}
			|(F_t^e)_{\bar z}(z)|
			&=\frac{|\HH(f_z(z+te))-\HH(f_z(z))|}{|t|}\\
			&\le k\frac{|f_z(z+te)-f_z(z)|}{|t|}
			=k|(F_t^e)_z(z)|
			\qquad\text{for a.e. }z\in U'.
		\end{align*}
		The $L^2$ bound needed to apply Caccioppoli follows directly from
		absolute continuity of Sobolev functions on almost every line parallel
		to $e$. The fundamental theorem of calculus on such lines gives
		\[
		F_t^e(z)=\int_0^1\partial_ef(z+ste)\,\dd s.
		\]
		Jensen's inequality and a translation of the area variable imply
		\[
		\int_{U'}|F_t^e|^2\dA
		\le\int_0^1\int_{U'}|Df(z+ste)|^2\,\dd A(z)\,\dd s
		\le\int_{U''}|Df|^2\dA.
		\]
		In particular, this bound is independent of $t$. Choose a real cutoff
		$\eta\in C_c^\infty(U')$ equal to $1$ on a neighborhood of
		$\overline U$. Lemma~\ref{lem:caccioppoli}, with $w=F_t^e$ and
		$q=0$, now yields
		\[
		\int_U|DF_t^e|^2\dA
		\le C_k\|\nabla\eta\|_\infty^2
		\int_{U'}|F_t^e|^2\dA
		\le C_{k,U,U'}\int_{U''}|Df|^2\dA.
		\]
		For fixed $t$, $F_t^e$ already belongs to $W^{1,2}(U')$ because
		$f\in W^{1,2}(U'')$. Moreover, $DF_t^e$ is the difference quotient
		of $Df$ in the direction $e$. Applying the difference-quotient
		criterion to $Df$ for $e=(1,0)$ and $e=(0,1)$ proves that
		$Df\in W^{1,2}(U)$. Since $U\Subset\Omega$ is arbitrary, this gives
		$f\in W^{2,2}_{\loc}$.
		This regularity reduction is established in
\cite[Proposition~2.1]{ACFJK17} and \cite[Proposition~8]{ACFJK20}.
		
		Set $g=f_z$ and $h=f_{\bar z}$. Both maps are now in
		$W^{1,2}_{\loc}$, and $h=\HH(g)$ almost everywhere. On almost every
		line parallel to a fixed $e$, choose their absolutely continuous
		representatives. They satisfy $h=\HH(g)$ along the line: the equality
		holds almost everywhere there and extends to every interior point by
		continuity. At points where the one-dimensional derivatives of $g$
		and $h$ exist, apply
		\[
		|h(z+se)-h(z)|\le k|g(z+se)-g(z)|
		\]
		and divide by $|s|$ before letting $s\to0$. This proves
		$|\partial_eh|\le k|\partial_eg|$ almost everywhere. No derivative of
		the Lipschitz field $\HH$ is required.
		
		For each fixed direction the exceptional null set may be different.
		To obtain one set on which all directional inequalities hold, take
		a countable dense subset $\{e_j\}_{j\ge1}$ of the unit circle and
		remove the union of its exceptional sets, together with the null set
		where the weak Cartesian derivatives of $g$ and $h$ are not defined.
		At every remaining point, the expressions
		$\partial_eg=e_xg_x+e_yg_y$ and
		$\partial_eh=e_xh_x+e_yh_y$ depend continuously on $e$. Approximating
		any unit $e$ by a sequence from $\{e_j\}$ proves
		\eqref{eq:directional} there. This is the common full-measure
		formulation of \cite[Proposition 8]{ACFJK20}.
	\end{proof}
	
	\subsection{A quadratic form for the second derivatives}
	For a map satisfying \eqref{eq:directional}, write
	\begin{equation}\label{eq:hg}
		g=f_z,\quad h=f_{\bar z},\qquad
		u=f_{zz},\quad v=f_{z\bar z},\quad w=f_{\bar z\bar z}.
	\end{equation}
	The next lemma gives both the sign of the quadratic form used in the
	proof and a way to estimate the full Hessian from its tangential
	components.
	
	\begin{lemma}\label{lem:calibration}
		Suppose that $f$ satisfies \eqref{eq:directional}, and use the notation
		\eqref{eq:hg}. Then
		\begin{equation}\label{eq:modulus}
			|w|+(1-k)|v|\le k|u|.
		\end{equation}
		For each $c\in\C$, it follows that
		\begin{equation}\label{eq:Q}
			Q_c(f):=|w-(1-k)e^{2i\arg(z-c)}v|^2-k^2|u|^2\le0
			\quad\text{almost everywhere}.
		\end{equation}
		Moreover, for every orthonormal frame $R,T$,
		\begin{equation}\label{eq:tangent-completion}
			|D^2f|^2\le2(K^2+2)
			\bigl(|\partial_Tg|^2+|\partial_Th|^2\bigr),
			\qquad K=\frac{1+k}{1-k}.
		\end{equation}
	\end{lemma}
	\begin{proof}
		By the formula for $\partial_e$, the hypothesis \eqref{eq:directional}
		reads
		\[
		|ve+w\bar e|\le k|ue+v\bar e|\qquad(|e|=1).
		\]
		The common full-measure set permits the choice of $e$ to depend on
		the point. Choose $e$ for which the left-hand side equals $|v|+|w|$.
		The right-hand side is at most $k(|u|+|v|)$. Rearranging gives
		\eqref{eq:modulus}, as in \cite[Lemma 11]{ACFJK20}.
		The triangle inequality then gives
		\eqref{eq:Q}. The phase can be assigned any value at $z=c$.
		
		Fix the frame at the point where the Hessian is evaluated. For the
		constant direction $R$, the map $F=\partial_Rf$ satisfies the
		quasiregular differential inequality. In the rotated coordinates, this
		reads
		\[
		|F_R+iF_T|\le k|F_R-iF_T|.
		\]
		The triangle inequality implies $|F_R|\le K|F_T|$, or
		$|f_{RR}|\le K|f_{RT}|$. The Frobenius norm is invariant under
		rotation, and
		\[
		|\partial_Tg|^2+|\partial_Th|^2
		=\tfrac12(|f_{RT}|^2+|f_{TT}|^2),
		\]
		so \eqref{eq:tangent-completion} follows. The derivatives of $F$ are
		taken in constant directions, and no derivative of the polar frame
		occurs.
	\end{proof}
	
	\subsection{Fourier series in angle and the Mellin variable}
	The weighted estimate will first be proved for an arbitrary complex-valued
	test map $\varphi\in C_c^\infty(\C,\C)$. No differential constraint
	is imposed on $\varphi$ in this calculation. We first place the center
	of the weight at zero. Translation then gives the formulas for a
	general center $c$. In this notation, the quadratic expression in \eqref{eq:Q}
	means
	\[
	Q_0(\varphi)(z)
	=\left|\varphi_{\bar z\bar z}(z)
	-(1-k)e^{2i\vartheta}\varphi_{z\bar z}(z)\right|^2
	-k^2|\varphi_{zz}(z)|^2,
	\qquad z=re^{i\vartheta}\ne0.
	\]
	The subscript $0$ specifies the center of the angular factor.
	The hypothesis on the solution $f$ gives $Q_0(f)\le0$, whereas $Q_0(\varphi)$ for
	an arbitrary test map need not have this sign. Its weighted integral
	will instead be estimated from below in Lemma~\ref{lem:test}.
	
	Fix $0<a<1$ and put $\mathcal{A}=1+a$. Before taking a Fourier transform,
	we remove the full real affine Taylor polynomial of $\varphi$ at zero:
	\[
	\ell(z)=\varphi(0)+\varphi_x(0)x+\varphi_y(0)y
	=\varphi(0)+\varphi_z(0)z+\varphi_{\bar z}(0)\bar z,
	\qquad \varphi^\circ=\varphi-\ell.
	\]
	Although the original test map $\varphi$ is compactly supported,
	$\varphi^\circ$ generally is not. The forms under consideration are
	unchanged because
	\[
	D^2\varphi^\circ=D^2\varphi,\qquad
	Q_0(\varphi^\circ)=Q_0(\varphi).
	\]
	The logarithmic radial coordinate and the weighted representative
	$V$ of $\varphi^\circ$ are defined by
	\begin{equation}\label{eq:mellin-scaling}
		\begin{aligned}
			&s=\log r,\qquad
			V(s,\vartheta)=e^{-(1+a)s}\varphi^\circ(e^{s+i\vartheta}),\\
			&\varphi^\circ(re^{i\vartheta})=r^{1+a}V(s,\vartheta).
		\end{aligned}
	\end{equation}
	Here $s\in\R$ and $\vartheta\in\T$: the punctured plane is written
	as the cylinder $\R\times\T$, with each circle $|z|=r$
	corresponding to one cross-section $s=\log r$. A dilation of $r$
	becomes a translation of $s$. The function $V$ is a rescaled
	representation of the test map, rather than a new solution of the
	Beltrami equation.
	
	The exponent $1+a$ in \eqref{eq:mellin-scaling} is chosen to
	remove the radial weight from the Hessian energy. Indeed,
	$\partial_r=r^{-1}\partial_s$, and each Cartesian derivative is
	$r^{-1}$ times an operator involving $\partial_s$, $\partial_\vartheta$
	and bounded angular coefficients. Applied twice to $r^{1+a}V$,
	these derivatives give $r^{a-1}$ times a differential expression
	in $V$ on the cylinder. Squaring gives $r^{2a-2}$, while
	\[
	\dd A=r\,\dd r\,\dd\vartheta
	=r^2\,\dd s\,\dd\vartheta,\qquad
	r^{-2a}r^{2a-2}r^2=1.
	\]
	Consequently, the weighted quadratic integral becomes an unweighted
	integral on the cylinder. We shall compute its differential
	coefficients explicitly below.
	
	The affine subtraction also justifies the transform at both ends of
	the cylinder. The map $\varphi^\circ$ is $O(r^2)$ at zero and grows
	at most linearly at infinity. Its representative $V$ and all its
	derivatives of order at most two decay like $e^{(1-a)s}$ as
	$s\to-\infty$ and at least like $e^{-as}$ as $s\to+\infty$.
	In particular, $V\in H^2(\R\times\T)$.
	Here $H^2(\R\times\T)$ consists of functions whose weak partial
	derivatives of total order at most two are square integrable with
	respect to $\dd s\,\dd\vartheta$.
	
	The connection with $f$ is as follows. In the subcritical proof we
	take a cutoff of the solution after subtracting a real affine map:
	\[
	\varphi=\chi(f-L),\qquad
	L(z)=f_{B(c,R)}+(Df)_{B(c,R)}(z-c).
	\]
	Here $\chi$ is a real smooth cutoff supported in $B(c,R)$ and
	equal to one on $B(c,R/2)$.
	The resulting map $\varphi$ belongs to $W^{2,2}$. Its use in the smooth test estimate
	is justified by the averaging over centers and approximation in
	Proposition~\ref{prop:subcritical}. In the critical argument, we
	instead subtract the affine Taylor polynomial of $f$ at $c$ and set
	\[
	U(z)=f(z)-f(c)-Df(c)(z-c),\qquad
	V(s,\vartheta)=e^{-(1+\gamma)s}U(c+e^{s+i\vartheta}).
	\]
	This is the same rescaling with $a=\gamma$, now applied locally to
	$U$. See \eqref{eq:Taylor}--\eqref{eq:critical-cylinder}. In that
	argument $V$ is initially in $H^2$ on finite subcylinders only, and
	its behavior as $s\to-\infty$ must be proved separately. The
	global $H^2$ observation above concerns the normalized smooth test
	map used in the present Fourier calculation.
	
	We first expand $V$ in the periodic variable $\vartheta$, and then
	take the one-dimensional Fourier transform of each coefficient in
	$s$. Our conventions are
	\[
	c_m(s)=\frac1{2\pi}\int_0^{2\pi}
	V(s,\vartheta)e^{-im\vartheta}\,\dd\vartheta,
	\qquad
	\widehat c_m(\xi)=\frac1{\sqrt{2\pi}}
	\int_\R c_m(s)e^{-i\xi s}\,\dd s,
	\qquad m\in\Z.
	\]
	The coefficient $c_m(s)$ gives the contribution of angular frequency
	$m$ on the circle of radius $e^s$. Its Fourier transform
	$\widehat c_m(\xi)$ resolves the dependence on logarithmic radius
	into frequencies $\xi\in\R$. The second integral is interpreted as an $L^2$
	Fourier transform when it is not absolutely convergent. Angular
	orthogonality and one-dimensional Plancherel give, in succession,
	\[
	\int_{\R\times\T}|V|^2\,\dd s\,\dd\vartheta
	=2\pi\sum_{m\in\Z}\int_\R|c_m(s)|^2\,\dd s
	=2\pi\sum_{m\in\Z}\int_\R|\widehat c_m(\xi)|^2\,\dd\xi.
	\]
	The factor $2\pi$ is $\int_0^{2\pi}|e^{im\vartheta}|^2\,\dd\vartheta$.
	The normalization $1/\sqrt{2\pi}$ makes the Fourier transform in
	$s$ preserve the $L^2$ norm, so no further constant occurs in the
	second equality. This identity also applies to the derivatives of
	$V$ that occur below.
	The transform sends $\partial_\vartheta$ to multiplication by $im$ and
	$\partial_s$ to multiplication by $i\xi$. Since
	$e^{i\xi s}=r^{i\xi}$, this is also a Mellin transform in $r$.
	The term Fourier--Mellin symbol refers here to the multiplier of a
	quadratic form after these two transforms. The calculation uses
	Fourier series and the ordinary one-dimensional Fourier transform.
	
	Put $\lambda=\mathcal{A}+i\xi$ and consider the mode
	\[
	\psi_{\lambda,m}(r,\vartheta)=r^\lambda e^{im\vartheta},\qquad
	r^\lambda=\exp(\lambda\log r),\quad r>0.
	\]
	Since $r^\lambda=r^{\mathcal{A}} e^{i\xi s}$, this is the original-plane
	expression corresponding to the Fourier mode
	$e^{i\xi s+im\vartheta}$ of $V$. The latter is not an
	$L^2(\R\times\T)$ function. It is used only to compute the
	differential multipliers. The integral identities below are applied
	to $V\in H^2(\R\times\T)$ through Plancherel's theorem.
	
	To compute the Hessian, the polar-coordinate chain rule first gives
	\[
	\partial_x=\cos\vartheta\,\partial_r
	-\frac{\sin\vartheta}{r}\partial_\vartheta,\qquad
	\partial_y=\sin\vartheta\,\partial_r
	+\frac{\cos\vartheta}{r}\partial_\vartheta.
	\]
	Consequently, the Wirtinger operators take the form
	\[
	\partial_z=\frac{e^{-i\vartheta}}2
	\left(\partial_r-\frac{i}{r}\partial_\vartheta\right),
	\qquad
	\partial_{\bar z}=\frac{e^{i\vartheta}}2
	\left(\partial_r+\frac{i}{r}\partial_\vartheta\right).
	\]
	Applying them once to $\psi=\psi_{\lambda,m}$ yields
	\[
	\begin{aligned}
		\psi_z
		&=\tfrac12(\lambda+m)r^{\lambda-1}e^{i(m-1)\vartheta},\\
		\psi_{\bar z}
		&=\tfrac12(\lambda-m)r^{\lambda-1}e^{i(m+1)\vartheta}.
	\end{aligned}
	\]
	Each differentiation lowers the radial exponent by one.
	Differentiation by $\partial_z$ also lowers the angular index by
	one, whereas $\partial_{\bar z}$ raises it by one. Thus, when
	$\partial_z$ is applied to $\psi_z$, its new coefficient is
	\[
	\frac{(\lambda-1)+(m-1)}2=\frac{\lambda+m-2}2.
	\]
	Similarly, applying $\partial_{\bar z}$ to $\psi_{\bar z}$
	gives the new coefficient
	$((\lambda-1)-(m+1))/2=(\lambda-m-2)/2$.
	For the mixed derivative, applying $\partial_z$ to
	$\psi_{\bar z}$ gives
	$((\lambda-1)+(m+1))/2=(\lambda+m)/2$.
	Multiplication by the respective first-derivative coefficients
	now gives all three entries of the Wirtinger Hessian:
	\begin{equation}\label{eq:wirtinger-symbols}
		\begin{aligned}
			\psi_{zz}
			&=\tfrac14(\lambda+m)(\lambda+m-2)
			r^{\lambda-2}e^{i(m-2)\vartheta},\\
			\psi_{z\bar z}
			&=\tfrac14(\lambda^2-m^2)r^{\lambda-2}e^{im\vartheta},\\
			\psi_{\bar z\bar z}
			&=\tfrac14(\lambda-m)(\lambda-m-2)
			r^{\lambda-2}e^{i(m+2)\vartheta}.
		\end{aligned}
	\end{equation}
	The factor $1/4$ comes from the two factors $1/2$ in the
	first-order operators. The terms $-2$ in the pure derivatives
	come from the changed radial exponent and angular index after
	the first differentiation. The mixed derivative retains angular
	index $m$ and has coefficient
	$(\lambda-m)(\lambda+m)/4=(\lambda^2-m^2)/4$.
	In the multiplier calculations we remove the common radial factor
	$r^{\lambda-2}$ from \eqref{eq:wirtinger-symbols}. The remaining
	angular factors are also needed: multiplication of
	$\psi_{z\bar z}$ by $e^{2i\vartheta}$ places it in the same angular
	mode as $\psi_{\bar z\bar z}$. This is why their difference in
	$Q_0$ can be represented by a single coefficient in
	Lemma~\ref{lem:symbol}.
	
	For comparison, \eqref{eq:polar-hessian} gives
	\[
	\begin{aligned}
		\psi_{RR}&=\lambda(\lambda-1)r^{\lambda-2}e^{im\vartheta},\\
		\psi_{RT}&=im(\lambda-1)r^{\lambda-2}e^{im\vartheta},\\
		\psi_{TT}&=(\lambda-m^2)r^{\lambda-2}e^{im\vartheta}.
	\end{aligned}
	\]
	These are evaluations of the Cartesian Hessian in the orthonormal
	frame $(R,T)$. They express the same second derivatives as
	\eqref{eq:wirtinger-symbols} in a different basis.
	
	We can now relate $Q_0$ directly to \eqref{eq:wirtinger-symbols}.
	For the mode $\psi=\psi_{\lambda,m}$ in that formula, put
	\[
	b_{m,a}(\xi)=\tfrac14(\lambda-m)
	[k\lambda-(2-k)m-2],\qquad
	d_{m,a}(\xi)=\tfrac{k}{4}(\lambda+m)(\lambda+m-2).
	\]
	The two expressions whose squared moduli define $Q_0$ are
	\[
	\begin{aligned}
		\psi_{\bar z\bar z}-(1-k)e^{2i\vartheta}\psi_{z\bar z}
		&=b_{m,a}(\xi)r^{\lambda-2}e^{i(m+2)\vartheta},\\
		k\psi_{zz}
		&=d_{m,a}(\xi)r^{\lambda-2}e^{i(m-2)\vartheta}.
	\end{aligned}
	\]
	In the first line, factoring $\lambda-m$ gives
	\[
	(\lambda-m-2)-(1-k)(\lambda+m)
	=k\lambda-(2-k)m-2,
	\]
	which explains its coefficient. Since $\Rea\lambda=1+a$,
	the pointwise value on this formal mode is
	\[
	Q_0(\psi)
	=r^{2a-2}\bigl(|b_{m,a}(\xi)|^2-|d_{m,a}(\xi)|^2\bigr).
	\]
	The radial factor disappears after multiplication by $r^{-2a}$
	and $\dd A=r^2\,\dd s\,\dd\vartheta$, as described above.
	Thus, the multiplier corresponding to the weighted quadratic form is
	\[
	P_{m,a}(\xi)=|b_{m,a}(\xi)|^2-|d_{m,a}(\xi)|^2.
	\]
	For a general $V$, the angular frequencies in the first line are
	$m+2$, and those in the second line are $m-2$. Different values of
	$m$ remain orthogonal within each squared modulus. Applying angular
	Parseval followed by Plancherel in $s$ gives
	\[
	\begin{aligned}
		\int_\C r^{-2a}Q_0(\varphi)\dA
		&=\int_\C r^{-2a}Q_0(\varphi^\circ)\dA\\
		&=2\pi\sum_{m\in\Z}\int_\R
		P_{m,a}(\xi)|\widehat c_m(\xi)|^2\,\dd\xi.
	\end{aligned}
	\]
	These identities hold for $V\in H^2$ by approximation and
	Plancherel's theorem \cite{SteinWeiss}. The individual plane modes
	serve only to compute the multipliers. The next lemma
	factorizes $P_{m,a}$. Define
	\begin{align}
		F_m(a)&=(1-k)(m^2-\mathcal{A})-a(1+k)m,\nonumber\\
		G_m(a)&=m^2-a(1-k)m-\mathcal{A}(1-ak),\label{eq:FGC}\\
		C_m(a)&=(1-k)(m+1)^2
		-ak\bigl((1+k)m+(1-k)\bigr).\nonumber
	\end{align}
	The corresponding tangential multiplier will contain
	\begin{equation}\label{eq:U}
		U_m(a)=(m-1)^2(\mathcal{A}+m)^2+(m+1)^2(\mathcal{A}-m)^2.
	\end{equation}
	
	\begin{lemma}\label{lem:symbol}
		For $m\in\Z$, $\xi\in\R$ and $\lambda=1+a+i\xi$, the multiplier
		defined above is
		\begin{align}
			P_{m,a}(\xi)
			&=\frac1{16}\left(
			|\lambda-m|^2|k\lambda-(2-k)m-2|^2
			-k^2|\lambda+m|^2|\lambda+m-2|^2\right)\nonumber\\
			&=\frac14\bigl(F_m(a)G_m(a)+C_m(a)\xi^2\bigr).
			\label{eq:symbol}
		\end{align}
		If $0<a<\gamma$, then $F_m(a)G_m(a)>0$ and $C_m(a)>0$ for every
		$m$. At $a=\gamma$, the product $F_m(\gamma)G_m(\gamma)$ vanishes
		only for $m=2$, and the coefficients $C_m(\gamma)$ remain positive.
		Hence, the only
		zero of $P_{m,\gamma}(\xi)$ is $(m,\xi)=(2,0)$.
	\end{lemma}
	\begin{proof}
		The first equality follows from the coefficients $b_{m,a}$ and
		$d_{m,a}$ computed above using \eqref{eq:wirtinger-symbols}.
		
		For the second equality, put
		\begin{align}
			r_m&=(\mathcal{A}-m)[k\mathcal{A}-(2-k)m-2],&
			\widetilde r_m&=k(\mathcal{A}+m)(\mathcal{A}+m-2),\nonumber\\
			\ell_m&=2k\mathcal{A}-2m-2,&
			\widetilde\ell_m&=2k(a+m).
			\label{eq:operator-coefficients}
		\end{align}
		The two polynomials inside the squared moduli are
		$r_m-k\xi^2+i\ell_m\xi$ and
		$\widetilde r_m-k\xi^2+i\widetilde\ell_m\xi$. Their difference of
		squared moduli has constant term $r_m^2-\widetilde r_m^2$ and
		coefficient of $\xi^2$ equal to
		$\ell_m^2-\widetilde\ell_m^2-2k(r_m-\widetilde r_m)$.
		The $\xi^4$ coefficients are equal and cancel. Expanding the
		coefficients gives
		\begin{equation}\label{eq:coefficient-identities}
			\begin{aligned}
				r_m-\widetilde r_m&=2F_m(a),&
				r_m+\widetilde r_m&=2G_m(a),\\
				\ell_m^2-\widetilde\ell_m^2-2k(r_m-\widetilde r_m)&=4C_m(a).
			\end{aligned}
		\end{equation}
		This proves the factorization.
		
		We now assume $0<a\le\gamma$. The four products that require a
		separate check are
		\begin{align*}
			F_{-1}G_{-1}&=2a^3k^2,\\
			F_0G_0&=(1-k)\mathcal{A}^2(1-ak),\\
			F_1G_1&=2a^2[2-k(a+2)],\\
			F_2G_2&=[3(1-k)-a(3+k)]
			[3(1-a)+ak(a+3)].
		\end{align*}
		The first three are positive: in the third one, use
		\[
		2-k(a+2)\ge2-k(\gamma+2)
		=\frac{(1-k)(6-k)}{3+k}>0.
		\]
		The fourth is positive for $a<\gamma$ and vanishes at $a=\gamma$.
		Its second factor is always positive in this range.
		
		Regarding $F_m$ and $G_m$ as polynomials of a real variable $m$,
		we have, for $m\ge2$,
		\[
		\frac{\dd F_m}{\dd m}
		=2(1-k)m-a(1+k)
		\ge\frac{(1-k)(9+k)}{3+k}>0,
		\qquad
		\frac{\dd G_m}{\dd m}=2m-a(1-k)>0.
		\]
		Thus, the signs at $m=2$ determine the signs for every integer
		$m\ge2$. On the half-line $m\le-1$, both derivatives are negative.
		Since $F_{-1}=2ak>0$ and $G_{-1}=a^2k>0$, both polynomials are
		positive throughout that half-line. This accounts for every angular
		mode.
		
		For $m\le-1$, we have
		\[
		C_m(a)\ge(1-k)(m+1)^2+2ak^2>0,
		\]
		and $C_0(a)=(1-k)(1-ak)>0$. For $m\ge1$, the inequalities
		$(m+1)^2\ge4m$ and $(1+k)m+(1-k)\le2m$ give
		\[
		C_m(a)\ge m[4(1-k)-2ak]>0.
		\]
		The last inequality follows from $ak\le k\gamma<2(1-k)$.
		In particular, $C_2(\gamma)>0$, so the vanishing constant term at
		$m=2$ gives a zero only when $\xi=0$.
	\end{proof}
	
	The zero $(2,0)$ corresponds to $V(s,\vartheta)=e^{2i\vartheta}$, or
	to the plane mode $r^{1+\gamma}e^{2i\vartheta}$. We call this the
	resonant mode. Its cylinder representative is constant in $s$ and
	does not belong to $L^2(\R\times\T)$. The vanishing multiplier describes
	the loss of a uniform bound at the critical weight. The flux argument
	in Section~\ref{subsec:flux} treats this mode on the half-cylinder.
	
	\begin{lemma}\label{lem:test}
		For $0<a<\gamma$ there is $\delta_{k,a}>0$ such that, for every
		$\varphi\in C_c^\infty(\C,\C)$ and every $c\in\C$,
		\begin{equation}\label{eq:test}
			\int_\C |z-c|^{-2a}Q_c(\varphi)\dA
			\ge\delta_{k,a}\int_\C |z-c|^{-2a}
			\bigl(|\partial_T\varphi_z|^2
			+|\partial_T\varphi_{\bar z}|^2\bigr)\dA.
		\end{equation}
	\end{lemma}
	\begin{proof}
		Translate $c$ to zero, keep $\varphi$ as the given compactly
		supported test map, and form $\varphi^\circ=\varphi-\ell$ as
		above. The Hessian and the tangential derivatives of
		$\varphi_z,\varphi_{\bar z}$ are unchanged: the first
		derivatives of $\ell$ are constants. Apply the Fourier identities
		to the $H^2$ representative $V$ of $\varphi^\circ$. On a mode
		$r^\lambda e^{im\vartheta}$,
		the angular derivatives of the scaled first Wirtinger derivatives
		have coefficients
		\[
		\tfrac i2(m-1)(\lambda+m),\qquad
		\tfrac i2(m+1)(\lambda-m).
		\]
		Since $\partial_T=r^{-1}\partial_\vartheta$ on each of
		$\varphi_z$ and $\varphi_{\bar z}$, the multiplier of the
		tangential integral in \eqref{eq:test} is
		\[
		\tfrac14\bigl(U_m(a)+2(m^2+1)\xi^2\bigr).
		\]
		For $a>0$, $U_m(a)$ is positive at every integer $m$. In addition,
		\[
		\frac{F_m(a)G_m(a)}{U_m(a)}\longrightarrow\frac{1-k}{2},
		\qquad
		\frac{C_m(a)}{2(m^2+1)}\longrightarrow\frac{1-k}{2}
		\quad\text{as }|m|\to\infty.
		\]
		Lemma~\ref{lem:symbol} gives
		\[
		\delta_{k,a}:=\inf_{m\in\Z}
		\min\left\{\frac{F_m(a)G_m(a)}{U_m(a)},
		\frac{C_m(a)}{2(m^2+1)}\right\}>0.
		\]
		Comparing the multipliers and applying Plancherel's theorem proves
		the result.
	\end{proof}
	
	Lemma~\ref{lem:test} estimates the tangential column of the Hessian.
	It does not estimate the full Hessian of an arbitrary test map: the
	full Hessian multiplier contains a term of order $\xi^4$, which has
	canceled in \eqref{eq:symbol}. In the next section, we first apply
	the test estimate to a cutoff and then use
	\eqref{eq:tangent-completion} where the cutoff agrees with the
	original solution.
	
	\section{Proofs of the main results}\label{sec:proofs}
	
	\subsection{Subcritical energy estimates}\label{subsec:subcritical}
	We first apply Lemma~\ref{lem:test} to a compactly supported cutoff of
	the solution. The singular weight need not be integrable against its
	Hessian at this stage. Averaging over the center gives bounded
	coefficients and allows us to use strong $W^{2,2}$ approximation.
	
	\begin{proposition}\label{prop:subcritical}
		Let $f$ satisfy the hypotheses of Theorem~\ref{thm:directional}.
		For $0<a<\gamma$ and $B(c,R)\Subset\Omega$, one has
		\begin{equation}\label{eq:subcritical-weight}
			\int_{B(c,R/2)}|z-c|^{-2a}|D^2f|^2\dA
			\le C_{k,a}R^{-2a}\int_{B(c,R)}|D^2f|^2\dA.
		\end{equation}
		Consequently, for $0<r\le R$,
		\begin{equation}\label{eq:subcritical-decay}
			\int_{B(c,r)}|D^2f|^2\dA
			\le C_{k,a}(r/R)^{2a}\int_{B(c,R)}|D^2f|^2\dA,
		\end{equation}
		and
		\begin{equation}\label{eq:subcritical-holder}
			[Df]_{a,B(c,R/2)}
			\le C_{k,a}R^{-a}
			\left(\int_{B(c,R)}|D^2f|^2\dA\right)^{1/2}.
		\end{equation}
	\end{proposition}
	\begin{proof}
		Choose a nonnegative smooth function $\rho_\varepsilon$, supported in
		$B(0,\varepsilon)$, with integral one. Integrate \eqref{eq:test} with
		respect to the centers $c+t$ against $\rho_\varepsilon(t)\dd A(t)$.
		The averaged weight
		\[
		W_\varepsilon(z)=\int_\C\rho_\varepsilon(t)
		|z-c-t|^{-2a}\dd A(t)
		\]
		is bounded for every fixed $\varepsilon>0$, because $a<1$ makes
		the kernel $|z|^{-2a}$ locally integrable. The coefficients of the averaged
		quadratic forms are bounded by $C_kW_\varepsilon$, since the angular
		phases and orthogonal projections have bounded modulus.
		
		The averaged inequality extends by strong $W^{2,2}$
		approximation to compactly supported Sobolev functions. To see this,
		view $H_j$ and $H$ as the arrays of second derivatives of the
		approximating maps and their limit, respectively. If
		$H_j\to H$ in $L^2(\C)$ and the averaged coefficient matrix
		$M_\varepsilon$ belongs to $L^\infty(\C)$, then
		\[
		\begin{aligned}
			&\left|\int_\C\langle M_\varepsilon H_j,H_j\rangle\dA
			-\int_\C\langle M_\varepsilon H,H\rangle\dA\right|\\
			&\quad\le\|M_\varepsilon\|_{L^\infty(\C)}
			\|H_j-H\|_{L^2(\C)}
			\bigl(\|H_j\|_{L^2(\C)}+\|H\|_{L^2(\C)}\bigr)
			\to 0.
		\end{aligned}
		\]
		
		Let $\chi\in C_c^\infty(B(c,R))$ equal one on $B(c,R/2)$, with
		$|D^j\chi|\le CR^{-j}$ for $j=1,2$, and put
		\[
		L(z)=f_{B(c,R)}+(Df)_{B(c,R)}(z-c),\qquad
		\varphi=\chi(f-L).
		\]
		The map $\varphi$, extended by zero outside $B(c,R)$, belongs to
		$W^{2,2}(\C,\C)$ and is admissible in the averaged inequality.
		The choice of $L$ gives
		$(f-L)_{B(c,R)}=0$ and $D(f-L)=Df-(Df)_{B(c,R)}$,
		since the mean of $z-c$ on the centered disk is zero. Poincar\'e's inequality applied
		first to $Df$ and then to $f-L$ gives
		\begin{align*}
			\int_{B(c,R)}|D(f-L)|^2\dA
			&\le CR^2\int_{B(c,R)}|D^2f|^2\dA,\\
			\int_{B(c,R)}|f-L|^2\dA
			&\le CR^4\int_{B(c,R)}|D^2f|^2\dA.
		\end{align*}
		In particular,
		\begin{equation}\label{eq:cutoff-bound}
			\|D^2\varphi\|_{L^2(\C)}
			\le C\|D^2f\|_{L^2(B(c,R))}.
		\end{equation}
		On $B(c,R/2)$, one has $D^2\varphi=D^2f$ and
		$Q_{c+t}(\varphi)\le0$. On the outer annulus
		$B(c,R)\setminus B(c,R/2)$, for
		$\varepsilon<R/8$,
		\[
		|z-c-t|\ge3R/8,\qquad
		|z-c-t|^{-2a}\le C_aR^{-2a},\qquad
		|Q_{c+t}(\varphi)|\le C_k|D^2\varphi|^2.
		\]
		Apply the averaged test inequality to $\varphi$ and retain the
		tangential energy on $B(c,R/2)$. The nonpositive contribution of
		$Q_{c+t}(\varphi)$ on this inner disk can be discarded from the upper
		bound. For each center $c+t$, \eqref{eq:tangent-completion} controls
		$|D^2f|^2$ by the retained energy, because $\varphi$ agrees with the
		solution there. Combining the outer-annulus bounds with
		\eqref{eq:cutoff-bound} gives
		\[
		\int_{B(c,R/2)}W_\varepsilon|D^2f|^2\dA
		\le C_{k,a}R^{-2a}\int_{B(c,R)}|D^2f|^2\dA.
		\]
		Since $W_\varepsilon(z)\to|z-c|^{-2a}$ for $z\ne c$, Fatou's
		lemma proves \eqref{eq:subcritical-weight}.
		For $r\le R/2$, the bound $|z-c|^{2a}\le r^{2a}$ on $B(c,r)$
		then gives \eqref{eq:subcritical-decay}. For $R/2<r\le R$, it
		follows from the total energy on $B(c,R)$.
		
		For every $y\in B(c,R/2)$, apply \eqref{eq:subcritical-decay} on
		$B(y,R/2)\subset B(c,R)$. Thus,
		\[
		\int_{B(y,\rho)}|D^2f|^2\dA
		\le C_{k,a}R^{-2a}\rho^{2a}\int_{B(c,R)}|D^2f|^2\dA,
		\qquad 0<\rho\le R/4.
		\]
		Lemma~\ref{lem:campanato}, applied to $Df$, proves
		\eqref{eq:subcritical-holder}.
	\end{proof}
	
	In what follows, $f$ denotes its $C^1$ representative. We shall use
	\eqref{eq:subcritical-holder} to estimate its Taylor remainder and
	\eqref{eq:subcritical-weight} to estimate the second derivatives of
	that remainder near the center.
	
	\subsection{A flux identity at the critical exponent}\label{subsec:flux}
	Fix a disk $B(c,R)\Subset\Omega$ and translate its center to zero.
	Subtract the affine Taylor polynomial:
	\begin{equation}\label{eq:Taylor}
		U(z)=f(z)-f(0)-Df(0)z.
	\end{equation}
	Here $Df(0)z$ is the real differential of $f$ acting on the vector
	$z\in\R^2$. In particular, $U(0)=0$ and $DU(0)=0$.
	Set
	\begin{equation}\label{eq:critical-cylinder}
		V(s,\vartheta)=e^{-(1+\gamma)s}U(e^{s+i\vartheta})
		=\sum_{m\in\Z}c_m(s)e^{im\vartheta},\qquad s<\log R,
	\end{equation}
	where
	\[
	c_m(s)=\frac1{2\pi}\int_0^{2\pi}
	V(s,\vartheta)e^{-im\vartheta}\dd\vartheta.
	\]
	The polar change of variables gives $V\in H^2$ on every finite
	subcylinder. We do not assume that $V$ belongs to $H^2$ on the whole
	half-cylinder. Its behavior as $s\to-\infty$ will instead be deduced
	from Proposition~\ref{prop:subcritical}.
	
	Evaluate the coefficients in \eqref{eq:operator-coefficients} at
	$a=\gamma$, and consequently write $\mathcal{A}=1+\gamma$. Define
	\begin{align}
		p_m&=r_m^2-\widetilde r_m^2
		=4F_m(\gamma)G_m(\gamma),\nonumber\\
		d_m&=\ell_m^2-\widetilde\ell_m^2-2k(r_m-\widetilde r_m)
		=4C_m(\gamma),\label{eq:pd}\\
		\Phi_m&=k(\ell_m-\widetilde\ell_m)|c_m'|^2
		+2k(r_m-\widetilde r_m)\Rea(c_m'\overline{c_m})\nonumber\\
		&\hspace{13mm}
		+(\ell_mr_m-\widetilde\ell_m\widetilde r_m)|c_m|^2.
		\label{eq:mode-flux}
	\end{align}
	We use the notation
	\begin{equation}\label{eq:bulk-flux}
		\DD(s)=\sum_{m\in\Z}(p_m|c_m|^2+d_m|c_m'|^2),\qquad
		\FF(s)=\sum_{m\in\Z}\Phi_m(s).
	\end{equation}
	The nonnegative expression $\DD$ is the energy density left after
	integration by parts in $s$. The quantity $\FF$ is the corresponding
	boundary term.
	At $m=2$, $p_m$ vanishes. The contribution of this mode to $\FF$ will
	provide the missing control of $c_2$.
	
	\begin{lemma}\label{lem:flux}
		The functions $\DD$ and $\FF$ belong to
		$L^1_{\loc}(( -\infty,\log R))$, and $\FF$ has a locally absolutely
		continuous representative. Define
		\begin{equation}\label{eq:q}
			q(s)=\frac1{2\pi}\int_0^{2\pi}
			e^{(2-2\gamma)s}Q_0(U)(e^{s+i\vartheta})\dd\vartheta\le0.
		\end{equation}
		Then
		\begin{equation}\label{eq:flux-identity}
			16q=\DD+\FF',\qquad \FF'\le-\DD\quad\text{a.e.}
		\end{equation}
		There is a constant $C_\Phi\ge1$, depending only on $k$, such that
		\begin{equation}\label{eq:flux-splitting}
			\FF=\BB-\MM,\qquad |\BB|\le C_\Phi\DD,
		\end{equation}
		where
		\begin{equation}\label{eq:resonant-mass}
			\MM=2(3+k)\bigl(k|c_2'|^2+R_*|c_2|^2\bigr),\qquad
			R_*=k(\gamma+3)(\gamma+1)>0.
		\end{equation}
		Moreover,
		\begin{equation}\label{eq:bulk-control}
			\sum_{m\ne2}(1+m^2)^2|c_m|^2
			+\sum_m(1+m^2)|c_m'|^2\le C_k\DD.
		\end{equation}
	\end{lemma}
	\begin{proof}
		Apply the cylinder change of variables to the two Hessian expressions
		defining $Q_0(U)$. In the $m$th angular mode, their scaled operators are
		\[
		\tfrac14(k\partial_s^2+\ell_m\partial_s+r_m),\qquad
		\tfrac14(k\partial_s^2+\widetilde\ell_m\partial_s+\widetilde r_m).
		\]
		The radial factor in each expression is $r^{\gamma-1}$ and is
		removed by the weight in \eqref{eq:q}. The angular factors have
		modulus one, so Parseval's identity gives
		\begin{equation}\label{eq:q-operators}
			16q=\sum_m\left(
			|kc_m''+\ell_mc_m'+r_mc_m|^2
			-|kc_m''+\widetilde\ell_mc_m'+\widetilde r_mc_m|^2\right).
		\end{equation}
		For a finite sum of angular modes, expand the squares and use
		\[
		\Rea(c_m''\overline{c_m'})=\tfrac12(|c_m'|^2)',\qquad
		\Rea(c_m''\overline{c_m})
		=(\Rea(c_m'\overline{c_m}))'-|c_m'|^2.
		\]
		The terms with $|c_m''|^2$ cancel. The remaining terms are
		$p_m|c_m|^2+d_m|c_m'|^2+\Phi_m'$. Summing proves the identity in
		\eqref{eq:flux-identity} for a finite Fourier series.
		
		By Lemma~\ref{lem:symbol}, $p_m>0$ for $m\ne2$, $p_2=0$, and
		$d_m>0$ for every $m$. Furthermore,
		\[
		p_m\sim4(1-k)m^4,\qquad d_m\sim4(1-k)m^2
		\quad\text{as }|m|\to\infty.
		\]
		These asymptotics, followed by comparison of the finitely many remaining
		coefficients, prove \eqref{eq:bulk-control}.
		For the exceptional mode, direct substitution gives
		\[
		r_2=\widetilde r_2=R_*,\qquad
		\ell_2-\widetilde\ell_2=-2(3+k).
		\]
		Hence, $\Phi_2=-\MM$. Define $\BB=\sum_{m\ne2}\Phi_m$.
		The coefficients of $|c_m'|^2$, $\Rea(c_m'\overline{c_m})$ and
		$|c_m|^2$ in \eqref{eq:mode-flux} are respectively
		$O(|m|)$, $O(m^2)$ and $O(|m|^3)$. Thus,
		\eqref{eq:bulk-control} also bounds $|\BB|$ by $C_\Phi\DD$.
		Equivalently, one may take $C_\Phi$ to be the maximum of one and the
		supremum, over $m\ne2$, of the operator norms of
		\[
		\begin{pmatrix}
			(\ell_mr_m-\widetilde\ell_m\widetilde r_m)/p_m
			& k(r_m-\widetilde r_m)/\sqrt{p_md_m}\\
			k(r_m-\widetilde r_m)/\sqrt{p_md_m}
			& k(\ell_m-\widetilde\ell_m)/d_m
		\end{pmatrix}.
		\]
		Each denominator is positive, and these matrices tend to zero as
		$|m|\to\infty$. Their supremum is finite.
		
		We now pass to the full Fourier series. Let $V_N$ be the sum over
		$|m|\le N$. On a finite subcylinder, $V_N\to V$ strongly in $H^2$.
		The two operator expressions converge in $L^2$, and their squared
		moduli converge in $L^1$. The preceding coefficient bounds also give
		$\DD_N\to\DD$ and $\FF_N\to\FF$ in $L^1$. For a finite
		sum, the product identities hold in distributions for $H^2$
		coefficients, by one-dimensional smooth approximation. Passing to the
		limit gives
		\[
		\FF'=16q-\DD\in L^1_{\loc}(( -\infty,\log R)).
		\]
		Thus, $\FF$ has a locally absolutely continuous representative.
		The sign $q\le0$ follows from \eqref{eq:Q}, applied to $U$, whose
		Hessian agrees with that of $f$. The Fourier truncations justify the
		identity. They are not required to satisfy the differential constraint,
		which is used only for the original solution to obtain $q\le0$ and
		the inequality in \eqref{eq:flux-identity}.
	\end{proof}
	
	\subsection{The flux at the inner end and the endpoint estimate}
	\label{subsec:endpoint}
	
	\begin{lemma}\label{lem:inner-sign}
		For the Taylor remainder \eqref{eq:Taylor}, the flux of
		Lemma~\ref{lem:flux} satisfies $\FF(s)\le0$ for every $s<\log R$.
		Consequently, $\WW=-\FF$ is nonnegative and nondecreasing, and
		\begin{equation}\label{eq:W-control}
			\WW(s)\le\WW(s_0)\quad(s\le s_0),\qquad
			\int_{-\infty}^{s_0}\DD(s)\dd s\le\WW(s_0).
		\end{equation}
	\end{lemma}
	\begin{proof}
		Put
		\begin{equation}\label{eq:epsilon-choice}
			c_k=(1+C_\Phi)^{-1},\qquad
			\varepsilon=\min\{\gamma/4,c_k/8\},\qquad
			a_1=\gamma-\varepsilon,\quad a_2=\gamma-\varepsilon/2.
		\end{equation}
		Both $a_1$ and $a_2$ belong to $(0,\gamma)$.
		We use $a_2$ to bound the Taylor remainder pointwise and $a_1$ to
		bound its weighted Hessian. The strict inequality $a_1<a_2$ leaves
		the integrable factor needed below.
		Proposition~\ref{prop:subcritical}, with exponent $a_2$, gives
		\[
		|U(z)|\le C|z|^{1+a_2},\qquad |DU(z)|\le C|z|^{a_2}
		\]
		near zero. The first bound follows by integrating the second along a
		segment, since $U(0)=DU(0)=0$. In cylinder coordinates,
		\begin{equation}\label{eq:low-order-inner}
			|V|+|V_s|+|V_\vartheta|\le C e^{-\varepsilon s/2}
			\qquad\text{for sufficiently negative }s.
		\end{equation}
		The polar Hessian formulas also give
		\begin{equation}\label{eq:cylinder-hessian}
			V_{\vartheta\vartheta}=r^{1-\gamma}U_{TT}-r^{-\gamma}U_R,
			\qquad
			V_{s\vartheta}=r^{1-\gamma}U_{RT}-\gamma r^{-\gamma}U_T.
		\end{equation}
		For a sufficiently small fixed $e^{s_1}$, the weighted estimate at
		exponent $a_1$ yields
		\[
		\int_{-\infty}^{s_1}\int_0^{2\pi}
		e^{2\varepsilon s}r^{2-2\gamma}|D^2U|^2\dd\vartheta\,\dd s
		=\int_{B(0,e^{s_1})}r^{-2a_1}|D^2U|^2\dA<\infty.
		\]
		After multiplication by $e^{2\varepsilon s}$, the squared lower-order
		terms in \eqref{eq:low-order-inner} and \eqref{eq:cylinder-hessian}
		are bounded by $Ce^{\varepsilon s}$. They are integrable at
		$-\infty$. From \eqref{eq:mode-flux},
		\begin{equation}\label{eq:flux-trace}
			|\FF(s)|\le C_k\bigl(\|V(s,\cdot)\|_{H^2(\T)}^2
			+\|V_s(s,\cdot)\|_{H^1(\T)}^2\bigr).
		\end{equation}
		It follows that
		\begin{equation}\label{eq:weighted-flux}
			\int_{-\infty}^{s_1}e^{2\varepsilon s}|\FF(s)|\dd s<\infty.
		\end{equation}
		Local integrability of $\FF$ extends this bound to every
		$s_1<\log R$. The constants may depend on the solution and the disk.
		Only the finiteness of the integral is needed to determine the sign
		of the flux.
		
		By \eqref{eq:flux-identity}, \eqref{eq:flux-splitting} and the choice
		$c_k C_\Phi<1$,
		\[
		\FF'+c_k\FF
		\le-\DD+c_k\BB-c_k\MM\le0.
		\]
		The function $e^{c_ks}\FF(s)$ is nonincreasing. If
		$\FF(s_*)>0$ at some $s_*<\log R$, then
		\[
		\FF(s)\ge e^{c_k(s_*-s)}\FF(s_*)\qquad(s<s_*).
		\]
		Multiplication by $e^{2\varepsilon s}$ gives a nonintegrable lower
		bound at $-\infty$, since $2\varepsilon<c_k$. This contradicts
		\eqref{eq:weighted-flux}. Thus, the absolutely continuous
		representative satisfies $\FF\le0$. Now
		$\WW'\ge\DD\ge0$, so for $t<s_0$,
		\[
		\int_t^{s_0}\DD(s)\dd s
		\le\WW(s_0)-\WW(t)\le\WW(s_0).
		\]
		Letting $t\to-\infty$ proves \eqref{eq:W-control}.
	\end{proof}
	
	\begin{lemma}\label{lem:annulus}
		For every $s_0<\log R$ and $s\le s_0$,
		\begin{equation}\label{eq:annular-critical}
			\int_{e^{s-1}<|z|<e^s}|z|^{-2\gamma}|D^2f|^2\dA
			\le C_k\WW(s_0).
		\end{equation}
	\end{lemma}
	\begin{proof}
		Let $I\subset(-\infty,s_0]$ be a finite interval of length at most one. Since
		$\MM=\BB+\WW$, \eqref{eq:flux-splitting} and
		\eqref{eq:W-control} imply
		\begin{equation}\label{eq:resonant-strip}
			\int_I\MM\dd s
			\le C_\Phi\int_I\DD\dd s+\int_I\WW\dd s
			\le(C_\Phi+1)\WW(s_0).
		\end{equation}
		Define
		\[
		G_\gamma V=\tfrac12e^{-i\vartheta}(\mathcal{A}+\partial_s-i\partial_\vartheta)V,
		\qquad
		H_\gamma V=\tfrac12e^{i\vartheta}(\mathcal{A}+\partial_s+i\partial_\vartheta)V.
		\]
		These expressions equal $r^{-\gamma}U_z$ and
		$r^{-\gamma}U_{\bar z}$. In the $m$th mode, their angular
		derivatives have coefficients
		\begin{equation}\label{eq:tangent-coefficients}
			\tfrac{i(m-1)}2\bigl(c_m'+(\mathcal{A}+m)c_m\bigr),\qquad
			\tfrac{i(m+1)}2\bigl(c_m'+(\mathcal{A}-m)c_m\bigr).
		\end{equation}
		For each $m$, the sum of the squared moduli of these coefficients is bounded by
		\[
		C_k\bigl((1+m^2)|c_m'|^2+(1+m^2)^2|c_m|^2\bigr).
		\]
		For $m\ne2$, its integral over $I$ is controlled by
		\eqref{eq:bulk-control} and \eqref{eq:W-control}. For $m=2$,
		\eqref{eq:resonant-mass} and \eqref{eq:resonant-strip} give the
		same bound. Parseval's identity yields
		\[
		\int_{I\times\T}
		\bigl(|\partial_\vartheta G_\gamma V|^2
		+|\partial_\vartheta H_\gamma V|^2\bigr)\dd\vartheta\,\dd s
		\le C_k\WW(s_0).
		\]
		Since $U$ differs from $f$ by an affine map, taking one more
		derivative removes that difference. With $g=f_z$ and $h=f_{\bar z}$,
		the two angular derivatives are
		$r^{1-\gamma}\partial_Tg$ and
		$r^{1-\gamma}\partial_Th$. Applying
		\eqref{eq:tangent-completion} to $f$ and taking $I=[s-1,s]$
		proves \eqref{eq:annular-critical}.
	\end{proof}
	
	\begin{proof}[Proof of Theorem~\ref{thm:directional}]
		We prove first that, whenever $B(c,R)\Subset\Omega$,
		\begin{equation}\label{eq:same-disk}
			\int_{B(c,r)}|D^2f|^2\dA
			\le C_k(r/R)^{2\gamma}\int_{B(c,R)}|D^2f|^2\dA,
			\qquad 0<r\le R.
		\end{equation}
		Translate $c$ to zero, let $U$ be the Taylor remainder in
		\eqref{eq:Taylor}, and put
		$E_R=\int_{B(0,R)}|D^2f|^2\dA$. With $a_0=\gamma/2$,
		\eqref{eq:subcritical-holder} gives
		\begin{equation}\label{eq:outer-Taylor}
			|DU(z)|\le C_k(r/R)^{a_0}E_R^{1/2},\qquad
			|U(z)|\le C_k r(r/R)^{a_0}E_R^{1/2},
			\qquad r=|z|<R/2.
		\end{equation}
		Again, integration along the radial segment gives the second bound.
		
		Consider
		\[
		J_R=[\log(R/4),\log(R/2)],\qquad |J_R|=\log2.
		\]
		On the corresponding annulus, the formulas \eqref{eq:cylinder-hessian} and
		the preceding Taylor estimates imply
		\begin{equation}\label{eq:outer-trace-average}
			\int_{J_R}\bigl(\|V(s,\cdot)\|_{H^2(\T)}^2
			+\|V_s(s,\cdot)\|_{H^1(\T)}^2\bigr)\dd s
			\le C_kR^{-2\gamma}E_R.
		\end{equation}
		To see the scaling, the second-derivative terms are bounded by
		\[
		C_k\int_{J_R\times\T}r^{2-2\gamma}|D^2f|^2
		\dd\vartheta\,\dd s
		\le C_kR^{-2\gamma}E_R.
		\]
		Every lower-order term is a multiple of $r^{-1-\gamma}U$ or
		$r^{-\gamma}DU$, and \eqref{eq:outer-Taylor} gives the same
		bound. Thus, \eqref{eq:outer-trace-average} uses the energy on an
		annulus rather than a trace on a fixed circle.
		
		Combining the averaged bound with \eqref{eq:flux-trace}, choose a point
		$s_0\in J_R$ where the trace bound holds and
		\begin{equation}\label{eq:outer-flux}
			\WW(s_0)\le C_kR^{-2\gamma}E_R.
		\end{equation}
		All constants depend only on $k$, and the same choice is available
		after translating to any center. For $r=e^s<R/4$, one has $s\le s_0$.
		Apply \eqref{eq:annular-critical} to the successive logarithmic
		intervals $[s-j-1,s-j]$, $j=0,1,\ldots$. Since
		$|z|^{2\gamma}\le e^{2\gamma(s-j)}$ on the $j$th annulus,
		\begin{align*}
			\int_{B(0,r)}|D^2f|^2\dA
			&\le C_k\WW(s_0)\sum_{j=0}^\infty e^{2\gamma(s-j)}\\
			&\le C_k r^{2\gamma}\WW(s_0)
			\le C_k(r/R)^{2\gamma}E_R.
		\end{align*}
		For $R/4\le r\le R$, the total energy $E_R$ gives the bound after
		increasing $C_k$. This proves \eqref{eq:same-disk}, which gives
		\eqref{eq:main-energy}.
		
		For $y\in B(c,R)$, apply \eqref{eq:same-disk} on
		$B(y,R)\subset B(c,2R)$. Then, for $0<\rho\le R/2$,
		\[
		\int_{B(y,\rho)}|D^2f|^2\dA
		\le C_kR^{-2\gamma}\rho^{2\gamma}
		\int_{B(c,2R)}|D^2f|^2\dA.
		\]
		Lemma~\ref{lem:campanato} yields \eqref{eq:main-holder}.
	\end{proof}
	
	\subsection{Natural solutions and the elliptic equation}
	\label{subsec:natural}
	
	\begin{proof}[Proof of Theorem~\ref{thm:natural} and sharpness]
		By Lemma~\ref{lem:natural}, a natural solution belongs to
		$W^{2,2}_{\loc}$ and satisfies \eqref{eq:directional}.
		Theorem~\ref{thm:directional} gives its regularity and the quantitative
		estimates.
		
		For sharpness, consider the example of
		\cite[Remark 15]{ACFJK20}, in the present normalization:
		\begin{equation}\label{eq:extremizer}
			f_\gamma(z)=z^2|z|^{\gamma-1},\quad f_\gamma(0)=0,\qquad
			\HH_\gamma(\zeta)=
			\begin{cases}
				-\dfrac{k}{3}\dfrac{\zeta^3}{|\zeta|^2},&\zeta\ne0,\\[2mm]
				0,&\zeta=0.
			\end{cases}
		\end{equation}
		For $z=re^{i\vartheta}\ne0$,
		\begin{equation}\label{eq:extremizer-derivatives}
			(f_\gamma)_z=\frac{3+\gamma}{2}r^\gamma e^{i\vartheta},\qquad
			(f_\gamma)_{\bar z}=-\frac{1-\gamma}{2}r^\gamma e^{3i\vartheta}.
		\end{equation}
		Since $(1-\gamma)/(3+\gamma)=k/3$, these derivatives satisfy
		$(f_\gamma)_{\bar z}=\HH_\gamma((f_\gamma)_z)$.
		At each nonzero point, the real differential of $\HH_\gamma$ has
		singular values $k/3$ and $k$ in the radial and tangential directions.
		Integration along
		a line segment, divided at zero when the segment passes through it,
		gives $\Lip(\HH_\gamma)\le k$. The tangential differential shows
		that equality holds.
		
		The derivatives in \eqref{eq:extremizer-derivatives} extend
		continuously to zero. Away from zero, the Hessian is homogeneous of
		degree $\gamma-1$, so its square is locally integrable. To identify
		these derivatives weakly across the origin, integrate by parts on a
		punctured disk. The inner boundary term for a derivative of $Df_\gamma$
		is bounded by $Cr^{1+\gamma}$ and tends to zero. Hence,
		$f_\gamma\in W^{2,2}_{\loc}$. For every $\alpha>\gamma$,
		\[
		\frac{|(f_\gamma)_z(r)-(f_\gamma)_z(0)|}{r^\alpha}
		=\frac{3+\gamma}{2}r^{\gamma-\alpha}\longrightarrow\infty.
		\]
		Thus, the exponent of Theorem~\ref{thm:natural} is optimal.
	\end{proof}
	
	The same example shows why the critical estimate was obtained on
	annuli. Indeed,
	\begin{equation}\label{eq:critical-divergence}
		\int_{B(0,1)}|z|^{-2\gamma}|D^2f_\gamma|^2\dA=\infty:
	\end{equation}
	the radial integral is a positive constant times $\int_0^1\dd r/r$.
	Its cylinder representation is $V=e^{2i\vartheta}$, with $\DD=0$ and
	a strictly negative constant flux. The estimate
	\eqref{eq:annular-critical} includes this example, although
	the corresponding weighted energy on a full disk is infinite.
	When $k=0$, subtraction of $\HH(0)\bar z$ reduces the autonomous
	equation to the holomorphic case, for which the limiting local exponent
	$\gamma=1$ follows directly.
	
	\begin{proof}[Proof of Corollary~\ref{cor:leray-lions}]
		We use the autonomous correspondence of
		\cite[Theorem 5]{ACFJK20}. Identify $p=(p_1,p_2)$ with
		$p_1+ip_2$, put $a=\mathbf{A}(p)$, and define
		$T(p)=p+\mathbf{A}(p)$. If $b=\mathbf{A}(q)$, squaring shows
		that \eqref{eq:strong-ellipticity} is equivalent to
		\begin{equation}\label{eq:Cayley-inequality}
			|(p-a)-(q-b)|\le k|(p+a)-(q+b)|,
			\qquad k=\frac{K-1}{K+1},
		\end{equation}
		because $K+K^{-1}=2(1+k^2)/(1-k^2)$.
		
		We check that $T$ is a bijection on $\R^2$.
		The ellipticity inequality and Cauchy--Schwarz give
		\[
		K^{-1}|p-q|\le|\mathbf{A}(p)-\mathbf{A}(q)|\le K|p-q|.
		\]
		Substituting the lower bound into \eqref{eq:strong-ellipticity} gives
		\[
		\langle p-q,\mathbf{A}(p)-\mathbf{A}(q)\rangle
		\ge K^{-1}|p-q|^2.
		\]
		It follows that $T$ is $(1+K)$-Lipschitz and strongly monotone with
		constant $\mu=1+K^{-1}$. For each $y\in\R^2$, the map
		\[
		p\longmapsto p-\tau(T(p)-y),\qquad
		\tau=\mu/(1+K)^2,
		\]
		has squared Lipschitz constant at most
		$1-\mu^2/(1+K)^2<1$. Its unique fixed point solves $T(p)=y$.
		
		On an interior disk $B$, the divergence-free field
		$\mathbf{A}(\nabla u)\in L^2(B,\R^2)$ has a weak stream function
		$v\in W^{1,2}(B)$ such that
		\[
		v_x=-\mathbf{A}_2(\nabla u),\qquad
		v_y=\mathbf{A}_1(\nabla u).
		\]
		For $f=u+iv$, with $p=u_x+iu_y$ and $a=\mathbf{A}(p)$,
		\[
		f_z=\frac{\overline{p+a}}2,\qquad
		f_{\bar z}=\frac{p-a}2.
		\]
		Define
		\begin{equation}\label{eq:Cayley-field}
			\HH(\zeta)=\tfrac12(I-\mathbf{A})T^{-1}(2\overline\zeta).
		\end{equation}
		The inequality \eqref{eq:Cayley-inequality} gives
		$\Lip(\HH)\le k$, and $\HH(0)=0$. The derivative identities
		above yield $f_{\bar z}=\HH(f_z)$ on $B$.
		Theorem~\ref{thm:natural} proves
		$u=\Rea f\in C^{1,\gamma}_{\loc}(B)$. Every point of $\Omega$
		lies in such a disk, which proves the local conclusion.
		
		To show optimality in the scalar class, use \eqref{eq:extremizer}.
		For each $p\in\C$, the equation
		\[
		p=\overline\zeta+\HH_\gamma(\zeta)
		\]
		has a unique solution, since
		$\zeta\mapsto\overline p-\overline{\HH_\gamma(\zeta)}$ is a
		contraction. Set
		$\mathbf{A}_\gamma(p)=\overline\zeta-\HH_\gamma(\zeta)$.
		The Lipschitz bound for $\HH_\gamma$ implies
		\eqref{eq:Cayley-inequality} for $\mathbf{A}_\gamma$, which gives
		\eqref{eq:strong-ellipticity}. Writing
		$f_\gamma=u_\gamma+iv_\gamma$, one obtains
		\[
		\mathbf{A}_\gamma(\nabla u_\gamma)
		=((v_\gamma)_y,-(v_\gamma)_x).
		\]
		Its divergence vanishes distributionally. On the other hand,
		$u_\gamma(r,\vartheta)=r^{1+\gamma}\cos(2\vartheta)$,
		$\nabla u_\gamma(0)=0$, and
		$(u_\gamma)_x(r,0)=(1+\gamma)r^\gamma$.
		Thus, $u_\gamma$ belongs to no $C^{1,\alpha}$ class near zero
		with $\alpha>\gamma$.
	\end{proof}
	
	\subsection{Growth of entire solutions}\label{subsec:entire}
	
	\begin{proof}[Proof of Theorem~\ref{thm:entire}]
		Suppose $L=\mathcal E_\gamma(f)<\infty$, and choose
		$R_j\to\infty$ with $R_j^{-2\gamma}E_f(R_j)\to L^2$.
		Fix $c\in\C$ and $r>0$, and put $S_j=(R_j-|c|)/2$.
		For large $j$, $B(c,2S_j)\subset B(0,R_j)$ and $r\le S_j$.
		The energy estimate \eqref{eq:main-energy} gives, as $j\to\infty$,
		\begin{equation}\label{eq:entire-morrey}
			\int_{B(c,r)}|D^2f|^2\dA\le C_kL^2r^{2\gamma},
		\end{equation}
		since $R_j/S_j\to2$. The whole-plane form of
		Lemma~\ref{lem:campanato} gives $[Df]_{\gamma,\C}\le C_kL$.
		
		Conversely, suppose $M=[Df]_{\gamma,\C}<\infty$.
		Each of $w=f_x,f_y$ satisfies the quasiregular differential inequality.
		Use a cutoff supported in $B(c,2r)$ and equal to one on $B(c,r)$.
		Lemma~\ref{lem:caccioppoli}, with $q=w(c)$, yields
		\[
		\int_{B(c,r)}|Dw|^2\dA
		\le C_kr^{-2}\int_{B(c,2r)}|w-w(c)|^2\dA
		\le C_kM^2r^{2\gamma}.
		\]
		Summing the estimates for $f_x$ and $f_y$ gives
		$\int_{B(c,r)}|D^2f|^2\dA\le C_kM^2r^{2\gamma}$. This proves the
		last inequality in \eqref{eq:global-equivalence}. The middle inequality
		follows from the ordering of the liminf and limsup. The two implications
		just proved also show that either $L$ and $M$ are both finite or both
		are infinite. Thus, \eqref{eq:global-equivalence} holds in the extended
		nonnegative reals.
		
		For an entire autonomous solution, set
		$F(z)=f(z)-\HH(0)\bar z$. Then
		\[
		|F_{\bar z}|=|\HH(f_z)-\HH(0)|\le k|f_z|=k|F_z|.
		\]
		Both $F$ and its directional derivatives satisfy the quasiregular
		differential inequality. Apply Lemma~\ref{lem:caccioppoli} first
		to $F_x,F_y$ with $q=0$, and then to $F$ with $q=F(0)$. This gives
		\begin{align}
			\int_{B(0,r)}|D^2f|^2\dA
			&\le C_kr^{-2}\int_{B(0,2r)}|DF|^2\dA\nonumber\\
			&\le C_kr^{-4}\int_{B(0,4r)}|F-F(0)|^2\dA.
			\label{eq:double-caccioppoli}
		\end{align}
		
		If $V_0=\mathcal V_\gamma(f)<\infty$, choose $R_j\to\infty$
		realizing the liminf in \eqref{eq:growth-value} and take $r=R_j/4$
		in \eqref{eq:double-caccioppoli}. Since
		\[
		\sup_{B(0,R_j)}|F-F(0)|
		\le\sup_{B(0,R_j)}|f|+|f(0)|+|\HH(0)|R_j,
		\]
		division by $R_j^{1+\gamma}$ makes the last two terms tend to zero.
		It follows that $\mathcal E_\gamma(f)\le C_kV_0$, so
		$[Df]_{\gamma,\C}\le C_kV_0$.
		For the converse estimate, integrate the differential along a segment:
		\[
		|f(z)-f(0)-Df(0)z|
		\le\frac{[Df]_{\gamma,\C}}{1+\gamma}|z|^{1+\gamma}.
		\]
		After division by $R^{1+\gamma}$, the affine terms vanish as
		$R\to\infty$. This proves the remaining inequality in
		\eqref{eq:value-equivalence}, including the extended-real case.
		
		If a growth quantity vanishes, the corresponding inequalities give
		$[Df]_{\gamma,\C}=0$, so $f$ is affine. To verify sharpness, the
		nonaffine example \eqref{eq:extremizer} satisfies
		\[
		E_{f_\gamma}(R)=E_{f_\gamma}(1)R^{2\gamma},\qquad
		\sup_{B(0,R)}|f_\gamma|=R^{1+\gamma}.
		\]
		Its gradient has a finite global $\gamma$-H\"older seminorm.
		To verify this, consider $h(re^{i\vartheta})=r^\gamma\Phi(\vartheta)$
		with smooth $\Phi$, and write $d=|z-w|$ and $\rho=\max\{|z|,|w|\}$.
		If $d\ge\rho/2$, then
		\mbox{$|h(z)-h(w)|\le C\rho^\gamma\le Cd^\gamma$}.
		If $d<\rho/2$, the segment joining the points stays at distance at
		least $\rho/2$ from zero. Integration of $|Dh|\le C\rho^{\gamma-1}$
		along that segment gives
		\mbox{$|h(z)-h(w)|\le C\rho^{\gamma-1}d\le Cd^\gamma$}.
		Apply this observation to \eqref{eq:extremizer-derivatives}. The
		growth powers are attained by a nonaffine solution.
	\end{proof}
	
	\section*{Statements and Declarations}

\subsection*{Author contributions}
Both authors contributed equally to the mathematical development and preparation of the manuscript and approved the final version.

\subsection*{Funding}
The second author was supported by NSF of China (No.~12271189, 12671096), NSF of Guangdong Province (Grant No. 2026A1515012333), STU Scientific Research Initiation Grant NTF25017T, HQU teaching reform project HQJGKT2411, and Fujian Alliance of Mathematics (Grant No.~2023SXLMMS07).

\subsection*{Competing interests}
The authors have no competing interests to disclose.

\subsection*{Data availability}
No external datasets were used.

\subsection*{Use of AI-assisted tools}
During the preparation of this manuscript, the authors used OpenAI’s ChatGPT to assist with language editing, stylistic refinement, and as a discussion aid for
checking the formulation and consistency of certain ideas. All mathematical statements, proofs,
computations, references, and conclusions were independently reviewed and verified by the
authors. The authors take full responsibility for the content of the manuscript.


\begin{thebibliography}{99}

\bibitem{Astala94}
K. Astala,
\href{https://doi.org/10.1007/BF02392568}{Area distortion of quasiconformal mappings},
\emph{Acta Math.} \textbf{173} (1994), no.~1, 37--60.

\bibitem{ACFJ}
K. Astala, A. Clop, D. Faraco and J. J\"a\"askel\"ainen,
\href{https://doi.org/10.1007/s11854-019-0059-x}{Manifolds of quasiconformal mappings and the nonlinear Beltrami equation},
\emph{J. Anal. Math.} \textbf{139} (2019), no.~1, 207--238.

\bibitem{ACFJK17}
K. Astala, A. Clop, D. Faraco, J. J\"a\"askel\"ainen and A. Koski,
\href{https://doi.org/10.1016/j.anihpc.2016.10.008}{Nonlinear Beltrami operators, Schauder estimates and bounds for the Jacobian},
\emph{Ann. Inst. H. Poincar\'e Anal. Non Lin\'eaire}
\textbf{34} (2017), no.~6, 1543--1559.

\bibitem{ACFJK20}
K. Astala, A. Clop, D. Faraco, J. J\"a\"askel\"ainen and A. Koski,
\href{https://doi.org/10.1016/j.matpur.2020.06.005}{Improved H\"older regularity for strongly elliptic PDEs},
\emph{J. Math. Pures Appl.} \textbf{140} (2020), 230--258.

\bibitem{ACFJS}
K. Astala, A. Clop, D. Faraco, J. J\"a\"askel\"ainen and L. Sz\'ekelyhidi Jr.,
\href{https://doi.org/10.1093/imrn/rnr178}{Uniqueness of normalized homeomorphic solutions to nonlinear Beltrami equations},
\emph{Int. Math. Res. Not. IMRN} \textbf{2012} (2012), no.~18, 4101--4119.

\bibitem{AFS}
K. Astala, D. Faraco and L. Sz\'ekelyhidi Jr.,
\href{https://doi.org/10.2422/2036-2145.2008.1.01}{Convex integration and the $L^p$ theory of elliptic equations},
\emph{Ann. Sc. Norm. Super. Pisa Cl. Sci. (5)} \textbf{7} (2008), no.~1, 1--50.

\bibitem{AIM}
K. Astala, T. Iwaniec and G. Martin,
\emph{\href{https://www.jstor.org/stable/j.ctt7sjhd}{Elliptic Partial Differential Equations and Quasiconformal Mappings
in the Plane}}, Princeton Mathematical Series, vol.~48,
Princeton University Press, Princeton, NJ, 2009.

\bibitem{BK}
A. Baernstein II and L. V. Kovalev,
\href{https://doi.org/10.2422/2036-2145.2005.2.04}{On H\"older regularity for elliptic equations of non-divergence type in the plane},
\emph{Ann. Sc. Norm. Super. Pisa Cl. Sci. (5)} \textbf{4} (2005), no.~2, 295--317.

\bibitem{Campanato}
S. Campanato,
\href{https://numdam.org/item/ASNSP_1963_3_17_1-2_175_0/}{Propriet\`a di h\"olderianit\`a di alcune classi di funzioni},
\emph{Ann. Scuola Norm. Sup. Pisa Cl. Sci. (3)} \textbf{17} (1963), no.~1--2, 175--188.

\bibitem{FinnSerrin}
R. Finn and J. Serrin,
\href{https://doi.org/10.2307/1993129}{On the H\"older continuity of quasi-conformal and elliptic mappings},
\emph{Trans. Amer. Math. Soc.} \textbf{89} (1958), no.~1, 1--15.

\bibitem{Giaquinta}
M. Giaquinta,
\emph{\href{https://api.pageplace.de/preview/DT0400.9781400881628_A26691462/preview-9781400881628_A26691462.pdf}{Multiple Integrals in the Calculus of Variations and Nonlinear
Elliptic Systems}}, Annals of Mathematics Studies, vol.~105,
Princeton University Press, Princeton, NJ, 1983.

\bibitem{HZ}
Z. Hou and J.-F. Zhu,
Improved H\"older regularity for elliptic equations of non-divergence type in the plane,
preprint, \href{https://arxiv.org/abs/2604.11536}{arXiv:2604.11536}, 2026.

\bibitem{IS}
T. Iwaniec and C. Sbordone,
\href{https://doi.org/10.1016/S0294-1449(00)00058-5}{Quasiharmonic fields},
\emph{Ann. Inst. H. Poincar\'e Anal. Non Lin\'eaire}
\textbf{18} (2001), no.~5, 519--572.

\bibitem{Kovalev}
L. V. Kovalev,
\href{https://doi.org/10.1112/jlms/jdm008}{Quasiconformal geometry of monotone mappings},
\emph{J. Lond. Math. Soc. (2)} \textbf{75} (2007), no.~2, 391--408.

\bibitem{LL}
J. Leray and J.-L. Lions,
\href{https://doi.org/10.24033/bsmf.1617}{Quelques r\'esultats de Vi\v{s}ik sur les probl\`emes elliptiques non lin\'eaires
par les m\'ethodes de Minty--Browder},
\emph{Bull. Soc. Math. France} \textbf{93} (1965), 97--107.

\bibitem{LN}
F. Leonetti and V. Nesi,
\href{https://doi.org/10.1016/S0021-7824(97)89947-3}{Quasiconformal solutions to certain first order systems and the proof of a
conjecture of G. W. Milton},
\emph{J. Math. Pures Appl.} \textbf{76} (1997), no.~2, 109--124.

\bibitem{Morrey}
C. B. Morrey Jr.,
\href{https://doi.org/10.2307/1989904}{On the solutions of quasi-linear elliptic partial differential equations},
\emph{Trans. Amer. Math. Soc.} \textbf{43} (1938), no.~1, 126--166.

\bibitem{SteinWeiss}
E. M. Stein and G. Weiss,
\emph{\href{https://www.degruyterbrill.com/document/doi/10.1515/9781400883899/html}{Introduction to Fourier Analysis on Euclidean Spaces}},
Princeton Mathematical Series, vol.~32,
Princeton University Press, Princeton, NJ, 1971.

\bibitem{Vekua}
I. N. Vekua,
\emph{\href{https://www.sciencedirect.com/book/monograph/9780080096933/generalized-analytic-functions}{Generalized Analytic Functions}}, International Series of Monographs on Pure and Applied Mathematics,
vol.~25, Pergamon Press, Oxford, 1962.

\end{thebibliography}
\end{document}